\documentclass[11pt]{article}
\usepackage[T1]{fontenc}
\usepackage{lmodern,amsmath,amssymb,amsthm,microtype}
\usepackage[a4paper,margin=29mm]{geometry}
\usepackage{graphicx,caption,placeins}
\usepackage[hidelinks]{hyperref}
\hypersetup{
 pdftitle={Uniqueness and nondegeneracy of ground states for the pure-power Born-Infeld equation},
 pdfkeywords={Born-Infeld equation, ground states, uniqueness, nondegeneracy, sharp Sobolev inequality}
}
\numberwithin{equation}{section}
\newtheorem{theorem}{Theorem}[section]
\newtheorem{lemma}[theorem]{Lemma}
\newtheorem{corollary}[theorem]{Corollary}
\newcommand{\R}{\mathbb R}
\newcommand{\D}{D^{1,2}}
\newcommand{\B}{\mathcal B}
\newcommand{\dx}{\,dx}
\newcommand{\dr}{\,dr}
\DeclareMathOperator{\diver}{div}
\DeclareMathOperator{\ind}{ind}
\title{Uniqueness and nondegeneracy of ground states\\
for the pure-power Born--Infeld equation}
\author{Huxiao Luo$\thanks{{\small Corresponding author. E-mail: luohuxiao@zjnu.edu.cn (H. Luo).}}$\\
\small\emph{School of Mathematical Sciences, Zhejiang Normal University, Jinhua, Zhejiang, {\rm321004}, China}
}
\date{}
\begin{document}
\maketitle
\begin{abstract}
For \(N\ge3\) and \(p>2N/(N-2)\), we study the Born--Infeld equation
\[
 -\diver\!\left(\frac{\nabla u}{\sqrt{1-|\nabla u|^2}}\right)
   =|u|^{p-2}u \qquad\text{in }\R^N .
\]
We prove that the positive ground state \(Q\) is unique up to
translation. We also prove that the kernel of the linearized operator
at \(Q\) is spanned by \(\partial_{x_1}Q,\ldots,\partial_{x_N}Q\),
and determine all nonzero optimizers of the sharp Born--Infeld Sobolev
inequality. To prove nondegeneracy, we first rule out nonzero radial
kernel functions using the Morse index and a Wronskian identity.
We then treat the angular modes by spherical harmonics. To prove
uniqueness, we compare radial solutions for large \(p\), and extend
the result to all \(p>2N/(N-2)\) by compactness and continuation.

\end{abstract}

\section{Introduction and main results}
Let \(N\ge3\), \(2^*=2N/(N-2)\), and \(p>2^*\). We consider
\begin{equation}\label{eq:pde}
 -\diver\!\left(\frac{\nabla u}{\sqrt{1-|\nabla u|^2}}\right)
 =|u|^{p-2}u \quad\text{in }\R^N .
\end{equation}
Let \(\D(\R^N)\) be the completion of \(C_c^\infty(\R^N)\)
under the norm \(\|\nabla u\|_2\). We identify this space with a
subspace of \(L^{2^*}(\R^N)\).
Following \cite{bim}, we call \(u\) admissible if
\[
 u\in\D(\R^N),\qquad \|\nabla u\|_\infty\le1 .
\]
Every admissible \(u\) has a unique \(1\)-Lipschitz representative
and belongs to \(L^\infty(\R^N)\).
We write \(\|\cdot\|_p\) for the \(L^p(\R^N)\) norm. Then
\[
 \|u\|_p^p\le \|u\|_\infty^{p-2^*}\|u\|_{2^*}^{2^*}<\infty
 \qquad(p>2^*).
\]
For each admissible \(u\), we define its energy by
\begin{equation}\label{eq:energy}
 I_p(u):=\int_{\R^N}\bigl(1-\sqrt{1-|\nabla u|^2}\bigr)\dx
 -\frac1p\int_{\R^N}|u|^p\dx .
\end{equation}
Both integrals are finite since \(0\le1-\sqrt{1-s^2}\le s^2\)
for \(0\le s\le1\) and \(u\in L^p(\R^N)\).

An admissible \(u\) is a finite-energy solution of \eqref{eq:pde} if
\begin{equation}\label{eq:weak}
 \int_{\R^N}\frac{\nabla u\cdot\nabla\varphi}
 {\sqrt{1-|\nabla u|^2}}\dx
 =\int_{\R^N}|u|^{p-2}u\varphi\dx
\end{equation}
holds for every \(\varphi\in C_c^\infty(\R^N)\), with both
integrals finite.
A ground state is a nonzero finite-energy solution of \eqref{eq:pde}
that minimizes \(I_p\) among all such solutions.

Bonheure, De Coster and Derlet \cite[Theorems 1.1 and 1.2]{bcd}
proved the existence of a positive radial finite-energy solution
for \(p>2^*\). They also proved that the best constant in a
Born--Infeld Sobolev inequality is attained by a radial solution.
Azzollini \cite[Theorems 0.2 and 0.3]{azz} found infinitely many
positive radial solutions with infinite Dirichlet energy and proved
the bound \(u(r)=O(r^{2-N})\) for positive radial solutions with
finite Dirichlet energy. In \cite{azz}, a ground state means a
positive solution that tends to zero at infinity. Throughout the
present paper, a ground state means a solution of least energy,
as defined above.

For other nonlinearities, Azzollini \cite{azz14} proved existence
by studying the radial initial value problem. Mederski and
Pomponio \cite{mp} obtained radial finite-energy solutions for
a class of nonlinearities. For a prescribed source, Bonheure and
Iacopetti \cite[Theorem 1.4]{bi} proved local \(W^{2,t}\) regularity,
with \(t>N\), for the energy minimizer under suitable assumptions
on the source. Byeon, Ikoma, Malchiodi and Mari \cite{bimm}
developed a variational method for autonomous Born--Infeld
equations and established regularity and Pohozaev identities.
Their Proposition~2.7 excludes nonzero solutions with finite
Dirichlet energy and finite \(L^p\) norm when \(2\le p\le2^*\).
Thus \(p>2^*\) is the full range of exponents \(p>2\) in which
finite-energy solutions of \eqref{eq:pde} can exist.

Bieganowski, Ikoma and Mederski \cite[Theorems 1.2 and 1.4]{bim}
proved the existence of ground states without imposing radial
symmetry in the minimization problem. For the pure-power equation,
they expressed the best Sobolev constant in terms of the
ground-state energy. They also proved that every nonzero optimizer
has the form \(sU(x/s)\), where \(s>0\) and \(U\) is a ground
state. These results describe optimizers in terms of ground states.
They do not prove uniqueness of positive ground states up to
translation or describe the kernel of the linearized operator.

We prove uniqueness up to translation and nondegeneracy of positive
ground states for every \(N\ge3\) and \(p>2^*\). The kernel
contains only the derivatives produced by translations.
Combining uniqueness with the optimizer characterization in
\cite{bim}, we describe all nonzero optimizers using a single
positive radial ground state, together with translations,
dilations, and a change of sign.

For a positive finite-energy solution \(Q\) of \eqref{eq:pde},
define the linearized operator by
\begin{equation}\label{eq:linear}
 L_Qh=-\diver\!\left(
 \frac{\nabla h}{\sqrt{1-|\nabla Q|^2}}
 +\frac{(\nabla Q\cdot\nabla h)\nabla Q}{(1-|\nabla Q|^2)^{3/2}}
 \right)-(p-1)Q^{p-2}h .
\end{equation}
Here \(h\in\ker_{\D(\R^N)}L_Q\) means that \(h\in\D(\R^N)\) and,
for every \(\varphi\in C_c^\infty(\R^N)\), we have
\[
 \int_{\R^N}\left(
 \frac{\nabla h\cdot\nabla\varphi}{\sqrt{1-|\nabla Q|^2}}
 +\frac{(\nabla Q\cdot\nabla h)(\nabla Q\cdot\nabla\varphi)}
 {(1-|\nabla Q|^2)^{3/2}}\right)\dx
 =(p-1)\int_{\R^N}Q^{p-2}h\varphi\dx .
\]
\begin{theorem}[Nondegeneracy]\label{th:nondegeneracy}
Let \(N\ge3\), \(p>2^*\), and let \(Q\) be any positive ground state of
\eqref{eq:pde}. Then
\begin{equation}\label{eq:kernel}
 \ker_{\D(\R^N)}L_Q
 =\operatorname{span}\{\partial_{x_1}Q,\ldots,\partial_{x_N}Q\}.
\end{equation}
\end{theorem}

\begin{theorem}[Uniqueness]\label{th:uniqueness}
Let \(N\ge3\) and \(p>2^*\). There is exactly one positive ground state
\(Q_p\) that is radial about the origin. The function \(Q_p\)
is strictly decreasing in
\(r=|x|\) for \(r>0\). Every positive ground state of \eqref{eq:pde}
has the form \(Q_p(\cdot-y)\) for some \(y\in\R^N\).
\end{theorem}

\begin{corollary}\label{cor:eq}
Let \(C_{N,p}\) be the optimal constant in the sharp Born--Infeld Sobolev inequality
\begin{equation}\label{eq:sharp}
 \int_{\R^N}\bigl(1-\sqrt{1-|\nabla u|^2}\bigr)\dx
 \ge C_{N,p}\,\|u\|_p^{Np/(N+p)} .
\end{equation}
A nonzero admissible function attains equality in
\eqref{eq:sharp} exactly when it has the form
\begin{equation}\label{eq:eqcases}
 u(x)=\sigma s\,Q_p\big((x-y)/s\big),
 \qquad \sigma\in\{-1,1\},\quad s>0,\quad y\in\R^N .
\end{equation}
\end{corollary}

The uniqueness criteria in Pucci and Serrin
\cite[Theorem 1]{ps} and Adachi, Shibata and Watanabe
\cite[Theorem 1.1]{asw} do not apply directly to \eqref{eq:pde}.
Both results require the nonlinearity to be negative for small
positive arguments, whereas \(f(s)=s^{p-1}>0\) for every \(s>0\).
The Born--Infeld operator is also not homogeneous: multiplying
the gradient by a constant changes the denominator in
\eqref{eq:pde}. The coefficient becomes unbounded as
\(|\nabla u|\) approaches one, so continuation requires uniform
gradient bounds for families of solutions. An amplitude and spatial rescaling cannot
normalize the value at the origin while preserving the equation.
The proof must control the dependence on the central value without
such a normalization.

For Theorem~\ref{th:nondegeneracy}, we first prove that positive
ground states are radial up to translation and have radial Morse
index one. We then suppose that a nonzero radial kernel function
exists. The index assumption forces the kernel function to have
exactly one zero. The equations for \(Q\) and the kernel function,
together with integration by parts, give two vanishing integrals
of the kernel function's radial flux, with weights \(q\) and
\(q^3\), where \(q=-Q'\). The two identities alone do not give
a contradiction, because the flux changes sign and \(q\) is not
monotone. A Wronskian identity shows that \(q\) has exactly one
maximum under the assumed existence of a radial kernel function.
The Wronskian identity also gives an increasing weight on two
intervals where \(q\) takes the same range of values. An exact
cancellation on these intervals lets us write a linear combination
of the two zero integrals as a sum of three strictly positive
terms. The resulting contradiction excludes radial kernel functions.
A spherical harmonic expansion then leaves only the translation
derivatives in the full kernel.

For Theorem~\ref{th:uniqueness}, we count all positive radial
finite-energy solutions with radial Morse index one. This class
contains the radial ground states. Nondegeneracy gives a unique
nearby solution when \(p\) changes slightly. Continuity of the
quadratic forms preserves the radial Morse index along each
local branch. Uniform bounds and compactness prevent additional
solutions from appearing outside these branches. The number of
solutions is therefore finite and constant on \((2^*,\infty)\).
Frank, Lenzmann and Silvestre \cite[Section 8]{fls} used a related
continuation argument in the order of the fractional Laplacian.
Here the varying parameter is the power \(p\), and the required
bounds must also control the nonlinear radial flux.

It remains to prove uniqueness for large \(p\). We first show
that the central height \(\alpha\) satisfies
\(cp^N\le\alpha^p\le Cp^N\), with constants independent of
the solution and of large \(p\). On this whole interval of
central values, differentiation of the radial initial value
problem gives the following ordering at a fixed radius \(r_0\):
increasing \(\alpha\) lowers the solution value and raises the
outward flux at \(r_0\). We also obtain a uniform bound
\(Q(r_0)\le\theta<1\). For large \(p\), this bound makes
the nonlinear term small enough on \([r_0,\infty)\) for an
energy estimate to rule out the ordering for two positive radial
finite-energy solutions. Uniqueness for large \(p\), together
with the constant number of solutions of radial Morse index one,
proves uniqueness of positive ground states for every \(p>2^*\).

\section{Proof of nondegeneracy}\label{sec:nond}
\subsection{Radial symmetry and the Morse index}\label{sec:symmetry-index}

We first give the regularity and energy identities of positive solutions.

\begin{lemma}\label{lem:regularity-energy}
Let \(Q\) be a positive finite-energy solution of \eqref{eq:pde}.
Then \(Q\in C^2(\R^N)\), \(\|\nabla Q\|_\infty<1\), and
\(Q(x)\to0\) as \(|x|\to\infty\).
Moreover,
\begin{equation}\label{eq:testq}
 \int_{\R^N}\frac{|\nabla Q|^2}{\sqrt{1-|\nabla Q|^2}}\dx=\|Q\|_p^p,
\end{equation}
and
\begin{equation}\label{eq:pohoz}
 I_p(Q)=\frac1N\|Q\|_p^p,\qquad
 \int_{\R^N}\bigl(1-\sqrt{1-|\nabla Q|^2}\bigr)\dx
 =\left(\frac1N+\frac1p\right)\|Q\|_p^p .
\end{equation}
\end{lemma}

\begin{proof}
We first check the weak formulation needed for the regularity
result. For any ball \(B\), choose \(\chi\in C_c^\infty(\R^N)\)
equal to one on \(B\), and test \eqref{eq:weak} with
\(\varphi(x)=x_i\chi(x)\). Finiteness of the integral gives
\(\partial_iQ/\sqrt{1-|\nabla Q|^2}\in L^1(B)\) for every
\(i=1,\ldots,N\). Since at least one component of \(\nabla Q\)
is nonzero where \(|\nabla Q|=1\), we obtain
\[
 |\nabla Q|<1\quad\text{a.e.},\qquad
 \frac{\nabla Q}{\sqrt{1-|\nabla Q|^2}}
 \in L^1_{\rm loc}(\R^N;\R^N).
\]
Also,
\[
 \frac1{\sqrt{1-|\nabla Q|^2}}
 \le 1+\frac{|\nabla Q|}{\sqrt{1-|\nabla Q|^2}}
 \in L^1_{\rm loc}(\R^N).
\]
Convolution approximates a compactly supported Lipschitz test
function by smooth test functions with uniformly bounded gradients
and a common compact support. The gradients converge almost
everywhere, so dominated convergence extends \eqref{eq:weak} to
compactly supported Lipschitz test functions. Thus \(Q\) satisfies
\cite[Definition~2.1]{bimm}.
Since \(Q\) and \(Q^{p-1}\) are bounded,
\cite[Proposition~2.10(ii)]{bimm} gives
\(\|\nabla Q\|_\infty<1\). Hence
\[
 \left|\frac{\nabla Q}{\sqrt{1-|\nabla Q|^2}}\right|
 \le C|\nabla Q|,
 \qquad \frac{\nabla Q}{\sqrt{1-|\nabla Q|^2}}\in L^2(\R^N).
\]
Since \(Q\in L^{2^*}\cap L^\infty\) and
\((p-1)2N/(N+2)>2^*\), we have
\(Q^{p-1}\in L^{2N/(N+2)}\).
By the Sobolev inequality and density, \eqref{eq:weak} extends
to every test function in \(\D(\R^N)\). Taking the test function
to be \(Q\) gives \eqref{eq:testq}.
By \eqref{eq:testq} and \cite[Proposition 2.6]{bimm},
\(Q\in W^{2,t}_{\rm loc}\) for \(2\le t<\infty\) and satisfies the Pohozaev identity
\begin{equation}\label{eq:pohoz-original}
 N I_p(Q)
 =\int_{\R^N}\frac{|\nabla Q|^2}{\sqrt{1-|\nabla Q|^2}}\dx .
\end{equation}
Combining \eqref{eq:testq} and \eqref{eq:pohoz-original} with the
definition of \(I_p\) gives \eqref{eq:pohoz}.
Taking \(t>N\), we find that the first derivatives of \(Q\)
are locally H\"older continuous. Since \(\|\nabla Q\|_\infty<1\),
the equation is uniformly elliptic and can be written as
\[
 -\sum_{i,j=1}^N\left(
 \frac{\delta_{ij}}{\sqrt{1-|\nabla Q|^2}}
 +\frac{\partial_iQ\,\partial_jQ}{(1-|\nabla Q|^2)^{3/2}}
 \right)\partial_{ij}Q=Q^{p-1}.
\]
The coefficients and the right-hand side are locally H\"older continuous,
so interior Schauder regularity gives \(Q\in C^2\).
Since \(Q>0\), the coefficients and the nonlinearity are smooth
on the range of \((Q,\nabla Q)\) in every compact set.
Repeated interior regularity therefore gives \(Q\in C^\infty_{\rm loc}\),
which allows the differentiations used later.
Since \(Q\) is Lipschitz and belongs to \(L^{2^*}(\R^N)\),
\(Q(x)\to0\) as \(|x|\to\infty\).
\end{proof}

For a positive radial finite-energy solution \(Q\) of
\eqref{eq:pde}, write \(r=|x|\) and \(q=-Q'\).
The radial form of \eqref{eq:pde} is
\begin{equation}\label{eq:radial-Q}
 -\frac{1}{r^{N-1}}
 \left(
 \frac{r^{N-1}Q'(r)}{\sqrt{1-Q'(r)^2}}
 \right)'
 =Q(r)^{p-1},
 \qquad r>0.
\end{equation}
To integrate \eqref{eq:radial-Q}, define the radial flux by
\begin{equation}\label{eq:radial-flux-def}
 m(r):=-\frac{r^{N-1}Q'(r)}{\sqrt{1-Q'(r)^2}}
       =\frac{r^{N-1}q(r)}{\sqrt{1-q(r)^2}},
 \qquad r>0,
\end{equation}
and extend \(m\) continuously to zero by \(m(0)=0\).
By the divergence theorem,
\[
 \int_{\partial B_r}\frac{-\nabla Q}{\sqrt{1-|\nabla Q|^2}}
       \cdot\nu\,dS
 =|\mathbb S^{N-1}|\,m(r)=\int_{B_r}Q^{p-1}\dx,
\]
where \(\nu\) is the outward unit normal. Thus \(m(r)\) is the
outward flux through \(\partial B_r\) divided by \(|\mathbb S^{N-1}|\).

Equations \eqref{eq:radial-Q} and \eqref{eq:radial-flux-def} give
\begin{equation}\label{eq:flux2}
\left\{
 \begin{aligned}
  m'(r)&=r^{N-1}Q(r)^{p-1},\\
 q(r)&=\frac{m(r)}{\sqrt{r^{2N-2}+m(r)^2}},
 \qquad r>0.
 \end{aligned} \right.
\end{equation}
Integrating the equation for \(m'\) gives
\begin{equation}\label{eq:flux1}
 m(r)=\int_0^r t^{N-1}Q(t)^{p-1}\,dt,
 \qquad r>0.
\end{equation}
Since \(Q>0\), \eqref{eq:flux1} gives \(m(r)>0\) for \(r>0\).
The definition \eqref{eq:radial-flux-def} then gives
\begin{equation}\label{eq:radial-monotonicity}
 Q'(r)<0,\qquad r>0.
\end{equation}
Applying \cite[Theorem 0.3]{azz} with power \(p-1\), we obtain
\begin{equation}\label{eq:Q-decay}
 Q(r)=O(r^{2-N})\qquad\text{as }r\to\infty.
\end{equation}
Since \((N-2)(p-1)>N\), the function \(r^{N-1}Q(r)^{p-1}\)
is integrable at infinity. Thus \eqref{eq:flux1} gives
\begin{equation}\label{eq:fast1}
 m(r)\longrightarrow
 M:=\int_0^\infty t^{N-1}Q(t)^{p-1}\,dt
 \in(0,\infty)
 \qquad\text{as }r\to\infty.
\end{equation}
Equation \eqref{eq:flux2} and differentiation with respect to \(r\) give
\begin{equation}\label{eq:fast2}
 \begin{aligned}
 q(r)&\le m(r)r^{1-N}=O(r^{1-N}),\\
 q'(r)&=(1-q(r)^2)^{3/2}
 \left(Q(r)^{p-1}-\frac{(N-1)m(r)}{r^N}\right)
 =O(r^{-N}),
 \end{aligned}
\end{equation}
as \(r\to\infty\).

We now prove radial symmetry of positive ground states.

\begin{lemma}\label{lem:radial-symmetry}
Let \(Q\) be a positive ground state of \eqref{eq:pde}.
There is \(y\in\R^N\) such that \(Q(\cdot+y)\) is radial
about the origin and strictly decreasing in \(r>0\).
\end{lemma}

\begin{proof}
By the energy identities
\eqref{eq:pohoz} and \cite[Theorem 1.4]{bim}, the sharp constant satisfies
\begin{equation}\label{eq:ground-equality}
 \frac{\displaystyle\int_{\R^N}
        \bigl(1-\sqrt{1-|\nabla Q|^2}\bigr)\dx}
      {\|Q\|_p^{Np/(N+p)}}
 =\left(\frac1N+\frac1p\right)\bigl(N I_p(Q)\bigr)^{p/(N+p)}=C_{N,p}.
\end{equation}
Let \(Q^*\) be the symmetric decreasing rearrangement of \(Q\). Then
\(\|Q^*\|_p=\|Q\|_p\). For \(t>0\), \((Q-t)_+\) has compact support
and rearrangement \((Q^*-t)_+\). By the P\'olya--Szeg\H{o} inequality
for \((Q-t)_+\) with exponent \(2j\), letting \(t\downarrow0\) gives
\begin{equation}\label{eq:rearr-norm}
 \|\nabla Q^*\|_{2j}\le\|\nabla Q\|_{2j}\qquad(j\ge1).
\end{equation}
By \eqref{eq:rearr-norm},
\[
 \|\nabla Q^*\|_{2j}
 \le\|\nabla Q\|_{\infty}^{1-1/j}\|\nabla Q\|_2^{1/j}.
\]
It follows that \(\|\nabla Q^*\|_\infty\le\|\nabla Q\|_\infty<1\).
For \(0\le s<1\), write
\[
 1-\sqrt{1-s^2}=\sum_{j\ge1}c_js^{2j},\qquad c_j>0.
\]
By monotone convergence, the rearrangement estimate \eqref{eq:rearr-norm},
and the sharp Born--Infeld Sobolev inequality \eqref{eq:sharp},
\begin{equation}\label{eq:rearr}
 \begin{aligned}
 C_{N,p}\|Q\|_p^{Np/(N+p)}
 &\le \int_{\R^N}\bigl(1-\sqrt{1-|\nabla Q^*|^2}\bigr)\dx \\
 &=\sum_{j\ge1}c_j\|\nabla Q^*\|_{2j}^{2j}
 \le\sum_{j\ge1}c_j\|\nabla Q\|_{2j}^{2j}\\
 &=\int_{\R^N}\bigl(1-\sqrt{1-|\nabla Q|^2}\bigr)\dx
 =C_{N,p}\|Q\|_p^{Np/(N+p)} .
 \end{aligned}
\end{equation}
Since \(c_j>0\), equality in \eqref{eq:rearr} forces equality
in \eqref{eq:rearr-norm} for every \(j\). In particular,
\begin{equation}\label{eq:20261106-e1}
\|\nabla Q^*\|_2=\|\nabla Q\|_2.
\end{equation}
Equation \eqref{eq:rearr} also shows that \(Q^*\) is an optimizer.
By \cite[Theorem 1.4]{bim}, \(s^{-1}Q^*(s\,\cdot)\) is a positive
ground state for some \(s>0\). The two ground states
\(s^{-1}Q^*(s\,\cdot)\) and \(Q\) have equal energy, so
\eqref{eq:pohoz} gives equal \(L^p\) norms. Hence
\[
 \|Q\|_p^p
 =\|s^{-1}Q^*(s\,\cdot)\|_p^p
 =s^{-(N+p)}\|Q^*\|_p^p
 =s^{-(N+p)}\|Q\|_p^p .
\]
Thus \(s=1\) and \(Q^*\) solves \eqref{eq:pde}.

Applying \eqref{eq:radial-monotonicity} to \(Q^*\) gives
\((Q^*)'(r)<0\) for \(r>0\).
For \(0<t<Q^*(0)\), write
\[
 \|\nabla Q\|_2^2=
 \|\nabla(Q-t)_+\|_2^2+\|\nabla\min(Q,t)\|_2^2 .
\]
The same identity holds for \(Q^*\). By the P\'olya--Szeg\H{o} inequality,
\[
 \begin{aligned}
 \|\nabla(Q^*-t)_+\|_2^2 &\le \|\nabla(Q-t)_+\|_2^2,\\
 \|\nabla\min\{Q^*,t\}\|_2^2 &\le \|\nabla\min\{Q,t\}\|_2^2.
 \end{aligned}
\]
Both inequalities are equalities, since otherwise their sum would
contradict \eqref{eq:20261106-e1}. In particular,
\[
 \|\nabla(Q-t)_+\|_2=\|\nabla(Q^*-t)_+\|_2 .
\]
Since \((Q^*)'(r)<0\) for \(r>0\), the gradient of \((Q^*-t)_+\) is
nonzero wherever its value is strictly between zero and its maximum.
The equality case of the rearrangement inequality
\cite{bz} (see also \cite[Theorem 1.1]{fv}) gives
\begin{equation}\label{eq:level-center}
 (Q-t)_+=(Q^*(\cdot-y_t)-t)_+
\end{equation}
for some \(y_t\). The equality holds everywhere by continuity.
For \(0<t_1<t_2<Q^*(0)\), applying \eqref{eq:level-center}
at \(t_1\) and \(t_2\) gives
\[
 \begin{aligned}
 \{Q>t_2\}
 &=\{(Q-t_1)_+>t_2-t_1\}
   =y_{t_1}+\{Q^*>t_2\},\\
 \{Q>t_2\}
 &=\{(Q-t_2)_+>0\}
   =y_{t_2}+\{Q^*>t_2\}.
 \end{aligned}
\]
Since \(\{Q^*>t_2\}\) is a nonempty bounded ball,
\(y_{t_1}=y_{t_2}\).
Thus all levels have one center \(y\), and \eqref{eq:level-center}
gives \(Q=Q^*(\cdot-y)\).
\end{proof}

For a positive radial finite-energy solution \(Q\), define
the bilinear form associated with \(L_Q\) by
\begin{equation}\label{eq:form}
 \B_Q(h,k)=\int_{\R^N}\left(
 \frac{\nabla h\cdot\nabla k}{\sqrt{1-|\nabla Q|^2}}
 +\frac{(\nabla Q\cdot\nabla h)(\nabla Q\cdot\nabla k)}
 {(1-|\nabla Q|^2)^{3/2}}-(p-1)Q^{p-2}hk\right)\dx .
\end{equation}
The coefficient matrix in \eqref{eq:form} is
\begin{equation}\label{eq:coefficient-matrix}
 \frac{I}{\sqrt{1-|\nabla Q|^2}}
 +\frac{\nabla Q\otimes\nabla Q}{(1-|\nabla Q|^2)^{3/2}} .
\end{equation}

Let \(\D_{\rm rad}\) denote the radial subspace of \(\D(\R^N)\).
The radial Morse index \(\ind(Q)\) is the largest dimension of a
subspace of \(\D_{\rm rad}\) on which \(\B_Q\) is negative definite.
We set \(\ind(Q)=+\infty\) if these dimensions are unbounded.

To prove that \(\ind(Q)\ge1\), we look for a nonzero radial function
\(Z\) such that \(\B_Q(Z,Z)<0\). Consider the scaled functions
\begin{equation}\label{eq:scaled-Q}
 Q_s(x):=sQ(x/s),\qquad s>0.
\end{equation}
We choose this scaling because \(\nabla Q_s(x)=\nabla Q(x/s)\).
We define \(Z\) as the derivative of \(Q_s\) at \(s=1\):
\begin{equation}\label{eq:Z-definition}
 \begin{aligned}
 Z(x)&:=\left.\frac{\partial Q_s(x)}{\partial s}\right|_{s=1}
       =Q(x)-x\cdot\nabla Q(x),\\
 Z(r)&=Q(r)-rQ'(r)=Q(r)+rq(r).
 \end{aligned}
\end{equation}
The next lemma shows that \(Z\) has the required property.

\begin{lemma}\label{lem:negative-direction}
Let \(Q\) be a positive radial finite-energy solution of \eqref{eq:pde},
and let \(Z\) be defined by \eqref{eq:Z-definition}.
Then \(\B_Q\) is continuous on \(\D(\R^N)\), \(Z\in\D_{\rm rad}\), and
\begin{equation}\label{eq:Z}
 Z>0,\qquad L_Q Z=-pQ^{p-1},\qquad
 \B_Q(Z,Z)=-p\int_{\R^N}Q^{p-1}Z\dx<0 .
\end{equation}
In particular, \(\ind(Q)\ge1\).
\end{lemma}

\begin{proof}
The matrix in \eqref{eq:coefficient-matrix} is bounded because \(\|\nabla Q\|_\infty<1\).
By \eqref{eq:Q-decay} and \((N-2)(p-2)>4\), we have
\(Q^{p-2}\in L^{N/2}\). H\"older's inequality gives
\[
 \left|\int_{\R^N}Q^{p-2}hk\dx\right|
 \le\|Q^{p-2}\|_{N/2}\|h\|_{2^*}\|k\|_{2^*}.
\]
The Sobolev inequality therefore proves continuity of \(\B_Q\)
on \(\D(\R^N)\).

For \(Q_s\) in \eqref{eq:scaled-Q}, equation \eqref{eq:pde} for \(Q\) gives
\begin{equation}\label{eq:scaled-equation}
 -\diver\!\left(
 \frac{\nabla Q_s(x)}{\sqrt{1-|\nabla Q_s(x)|^2}}
 \right)-Q_s(x)^{p-1}
   =(s^{-1}-s^{p-1})Q(x/s)^{p-1}.
\end{equation}
Differentiating \eqref{eq:scaled-equation} at \(s=1\) and using
\eqref{eq:Z-definition}, we obtain
\begin{equation}\label{eq1:2026118}
L_QZ=-pQ^{p-1}.
\end{equation}
The function \(Z\) is regular at zero. Definition
\eqref{eq:Z-definition} and the bounds \eqref{eq:Q-decay} and
\eqref{eq:fast2} give, as \(r\to\infty\),
\[
 Z=O(r^{2-N}),\qquad Z'=rq'=O(r^{1-N}).
\]
Thus \(Z\in\D_{\rm rad}\). Continuity of \(\B_Q\) and
\(Q^{p-1}\in L^{2N/(N+2)}\) allow us to test
\eqref{eq1:2026118} with \(Z\). Since \(Z=Q+rq>0\),
\[
 \B_Q(Z,Z)=-p\int_{\R^N}Q^{p-1}Z\dx<0.
\]
Hence \(\ind(Q)\ge1\).
\end{proof}

We now show that a radial ground state has Morse index one.

\begin{lemma}\label{lem:ground-morse}
Let \(Q\) be a positive ground state that is radial about the origin.
Then \(\ind(Q)=1\). Moreover, \(\B_Q(h,h)\ge0\) for every
\(h\in\D(\R^N)\) with \(\int_{\R^N}Q^{p-1}h\dx=0\).
\end{lemma}

\begin{proof}
By the sharp Born--Infeld Sobolev inequality \eqref{eq:sharp}
and equality at \(Q\) in \eqref{eq:ground-equality}, every
admissible \(v\) with \(\|v\|_p=\|Q\|_p\) satisfies
\begin{equation}\label{20261107-e1}
 \int_{\R^N}\bigl(1-\sqrt{1-|\nabla v|^2}\bigr)\dx
 \ge
 \int_{\R^N}\bigl(1-\sqrt{1-|\nabla Q|^2}\bigr)\dx.
\end{equation}

Let \(h\in C_c^\infty(\R^N)\) satisfy
\(\int_{\R^N}Q^{p-1}h\dx=0\). Define
\[
 a(t)=\frac{\|Q\|_p}{\|Q+th\|_p},
 \qquad
 v(t)=a(t)(Q+th).
\]
Then \(\|v(t)\|_p=\|Q\|_p\). Since \(p>2\), the function
\(a\) is \(C^2\) near zero. We have \(a(0)=1\), and the
condition on \(h\) gives \(a'(0)=0\). Hence
\[
 v(0)=Q,\qquad v'(0)=h.
\]
For sufficiently small \(|t|\),
\begin{equation}\label{eq:morse-path-gradient}
 \|\nabla v(t)\|_\infty
 \le |a(t)|\bigl(\|\nabla Q\|_\infty
                  +|t|\|\nabla h\|_\infty\bigr)
 \le \rho:=\frac{1+\|\nabla Q\|_\infty}{2}<1.
\end{equation}
Thus \(v(t)\) is admissible for small \(|t|\).

Differentiating the definition of \(v(t)\), we obtain
\[
 \begin{aligned}
 \partial_t v(t)
 &=a'(t)(Q+th)+a(t)h,\\
 \partial_t^2v(t)
 &=a''(t)(Q+th)+2a'(t)h.
 \end{aligned}
\]
Since \(a,a',a''\) are bounded near zero,
\begin{equation}\label{eq:morse-path-derivatives}
 \begin{aligned}
 &|\nabla v(t,x)|+|\nabla\partial_t v(t,x)|
   +|\nabla\partial_t^2v(t,x)|\\
 &\qquad\le C\bigl(|\nabla Q(x)|+|\nabla h(x)|\bigr).
 \end{aligned}
\end{equation}
Here and below, \(C\) is independent of \(t,x\).

Define
\[
 \Phi(z)=1-\sqrt{1-|z|^2},\qquad |z|<1,
 \qquad
 g(t)=\int_{\R^N}\Phi(\nabla v(t,x))\,dx.
\]
The matrix in \eqref{eq:coefficient-matrix} is
\(D^2\Phi(\nabla Q)\). For a general vector \(z\) with \(|z|<1\),
\[
 D\Phi(z)=\frac{z}{\sqrt{1-|z|^2}},
 \qquad
 D^2\Phi(z)=\frac{I}{\sqrt{1-|z|^2}}
       +\frac{z\otimes z}{(1-|z|^2)^{3/2}}.
\]
The formulas for \(D\Phi,D^2\Phi\) and
\eqref{eq:morse-path-gradient} give
\[
 |D\Phi(\nabla v)|\le C|\nabla v|,
 \qquad
 \|D^2\Phi(\nabla v)\|\le C.
\]
The chain rule gives
\[
 \begin{aligned}
 \partial_t \Phi(\nabla v)
 &=D\Phi(\nabla v)\cdot\nabla\partial_t v,\\
 \partial_t^2 \Phi(\nabla v)
 &=\nabla\partial_t v\cdot
   D^2\Phi(\nabla v)\nabla\partial_t v
   +D\Phi(\nabla v)\cdot\nabla\partial_t^2v.
 \end{aligned}
\]
The chain-rule formulas and \eqref{eq:morse-path-derivatives}
give
\[
 |\partial_t \Phi(\nabla v(t,x))|
 +|\partial_t^2 \Phi(\nabla v(t,x))|
 \le C\bigl(|\nabla Q(x)|+|\nabla h(x)|\bigr)^2.
\]
The right-hand side belongs to \(L^1(\R^N)\), and both
derivatives on the left are continuous in \(t\).
Dominated convergence therefore allows us to differentiate
\(g\) twice under the integral sign and gives \(g\in C^2\)
near zero.

Since \(v(t)\) is admissible and \(\|v(t)\|_p=\|Q\|_p\),
\eqref{20261107-e1} gives \(g(t)\ge g(0)\). Thus \(g''(0)\ge0\).
Using \(v(0)=Q\) and \(v'(0)=h\), we obtain
\[
 \begin{aligned}
 g''(0)
 &=\int_{\R^N}
 \left(
 \frac{|\nabla h|^2}{\sqrt{1-|\nabla Q|^2}}
 +\frac{(\nabla Q\cdot\nabla h)^2}
       {(1-|\nabla Q|^2)^{3/2}}
 \right)\dx\\
 &\quad+\int_{\R^N}
 \frac{\nabla Q\cdot\nabla v''(0)}
      {\sqrt{1-|\nabla Q|^2}}\dx.
 \end{aligned}
\]
To compute the integral containing \(v''(0)\), differentiate
\(\int_{\R^N}|v(t)|^p\dx=\|Q\|_p^p\) twice at \(t=0\).
Using \(v(0)=Q\) and \(v'(0)=h\), we obtain
\begin{equation}\label{eq:constraint-second}
 \int_{\R^N}Q^{p-1}v''(0)\dx
 =-(p-1)\int_{\R^N}Q^{p-2}h^2\dx.
\end{equation}
Since \(a'(0)=0\), the definition of \(v(t)\) gives
\(v''(0)=a''(0)Q\in\D(\R^N)\).
Testing \eqref{eq:weak} with \(v''(0)\) and using
\eqref{eq:constraint-second}, we obtain
\[
 \begin{aligned}
 \int_{\R^N}
 \frac{\nabla Q\cdot\nabla v''(0)}
      {\sqrt{1-|\nabla Q|^2}}\dx
 &=\int_{\R^N}Q^{p-1}v''(0)\dx\\
 &=-(p-1)\int_{\R^N}Q^{p-2}h^2\dx.
 \end{aligned}
\]
Consequently,
\begin{equation}\label{eq:tangent}
 \begin{aligned}
 0\le g''(0)
 &=\int_{\R^N}
 \left(
 \frac{|\nabla h|^2}{\sqrt{1-|\nabla Q|^2}}
 +\frac{(\nabla Q\cdot\nabla h)^2}
       {(1-|\nabla Q|^2)^{3/2}}
 -(p-1)Q^{p-2}h^2
 \right)\dx\\
 &=\B_Q(h,h).
 \end{aligned}
\end{equation}

We next extend \eqref{eq:tangent} to every \(h\in\D(\R^N)\)
with \(\int_{\R^N}Q^{p-1}h\dx=0\).
Choose \(h_i\in C_c^\infty(\R^N)\) such that \(h_i\to h\)
in \(\D\), and choose \(\varphi\in C_c^\infty(\R^N)\) with
\(\int_{\R^N}Q^{p-1}\varphi\dx=1\). Set
\[
 \widetilde h_i
 =h_i-\left(\int_{\R^N}Q^{p-1}h_i\dx\right)\varphi.
\]
Then \(\widetilde h_i\in C_c^\infty(\R^N)\) and
\[
 \int_{\R^N}Q^{p-1}\widetilde h_i\dx=0.
\]
Since \(Q^{p-1}\in L^{2N/(N+2)}(\R^N)\), the Sobolev
inequality and \(h_i\to h\) in \(\D\) give
\[
 \int_{\R^N}Q^{p-1}h_i\dx
 \longrightarrow
 \int_{\R^N}Q^{p-1}h\dx=0.
\]
Thus \(\widetilde h_i\to h\) in \(\D\).
Applying \eqref{eq:tangent} to \(\widetilde h_i\) and using
the continuity of \(\B_Q\), we obtain
\[
 \B_Q(h,h)\ge0
 \qquad\text{for every }h\in\D(\R^N)
 \text{ with }\int_{\R^N}Q^{p-1}h\dx=0.
\]

Every two-dimensional subspace of \(\D_{\rm rad}\) contains
a nonzero \(h\) satisfying \(\int_{\R^N}Q^{p-1}h\dx=0\).
The nonnegativity of \(\B_Q(h,h)\) rules out a two-dimensional
negative definite subspace. Thus \(\ind(Q)\le1\).
Lemma~\ref{lem:negative-direction} gives \(\ind(Q)\ge1\),
so \(\ind(Q)=1\).
\end{proof}

\subsection{The radial linearized equation}
We now rule out nonzero radial kernel functions when the radial
Morse index is one.

\begin{lemma}\label{lem:rad}
Let \(Q\) be a positive radial finite-energy solution of
\eqref{eq:pde} with \(\ind(Q)=1\). If a radial function \(\eta\in\D(\R^N)\)
satisfies \(L_Q\eta=0\), then \(\eta=0\).
\end{lemma}
\begin{proof}
We argue by contradiction. We first find the signs of a radial
kernel function \(\eta\) and its flux \(j\). We then show that
\(q=-Q'\) first increases and then decreases. The signs and
monotonicity imply that a linear combination of two zero integrals
is strictly positive. The figures illustrate the signs and turning
points used in the proof.

Suppose that \(0\ne\eta\in\D_{\rm rad}\) satisfies
\(L_Q\eta=0\). Recall \(q=-Q'\) and the flux \(m\) defined in
\eqref{eq:radial-flux-def}. The flux formulas
\eqref{eq:flux1} and \eqref{eq:flux2} give \(m(r)>0\) and
\(0<q(r)<1\) for \(r>0\).

\medskip\noindent
Step 1. We normalize \(\eta(0)=1\) and prove that \(\eta\) has
exactly one zero in \((0,\infty)\).

The radial form of \(L_Q\eta=0\) is
\begin{equation}\label{eq:radial-kernel}
 -\frac{1}{r^{N-1}}
 \left(r^{N-1}(1-q(r)^2)^{-3/2}\eta'(r)\right)'
 =(p-1)Q(r)^{p-2}\eta(r),\qquad r>0.
\end{equation}
For radial functions, the two gradient terms in
\eqref{eq:linear} have combined coefficient \((1-q^2)^{-3/2}\).
Lemma~\ref{lem:regularity-energy}, interior regularity for
\(\eta\), and radial symmetry give \(Q'(0)=\eta'(0)=0\). Multiplying \eqref{eq:radial-kernel}
by \(r^{N-1}\) and integrating from \(0\) to \(r\), with
zero boundary term at the origin, gives
\[
 -r^{N-1}(1-q(r)^2)^{-3/2}\eta'(r)
 =(p-1)\int_0^r t^{N-1}Q(t)^{p-2}\eta(t)\,dt.
\]
We denote the quantity on the left by \(j\):
\begin{equation}\label{eq:fluxlin}
 j(r)=-r^{N-1}(1-q(r)^2)^{-3/2}\eta'(r).
\end{equation}
The total outward flux of the linearized diffusion term through
\(\partial B_r\) is \(|\mathbb S^{N-1}|j(r)\).
The flux \(j\) lets us write \eqref{eq:radial-kernel} as
\begin{equation}\label{eq:j}
 \begin{aligned}
 j(0)&=0,\qquad
 j'(r)=(p-1)r^{N-1}Q(r)^{p-2}\eta(r),\\
 \eta'(r)&=-\frac{(1-q(r)^2)^{3/2}}{r^{N-1}}j(r)
 \qquad(r>0).
 \end{aligned}
\end{equation}
Thus \(j'\) has the same sign as \(\eta\), while \(\eta'\)
has the opposite sign to \(j\).

Integrating the equation for \(j'\) in \eqref{eq:j} and
substituting into the equation for \(\eta'\), we obtain, for small \(r>0\),
\[
 |j(t)|\le Ct^N\sup_{0\le s\le r}|\eta(s)|,
 \qquad
 |\eta'(t)|\le Ct\sup_{0\le s\le r}|\eta(s)|
 \qquad(0<t\le r).
\]
If \(\eta(0)=0\), integration gives
\[
 \sup_{0\le t\le r}|\eta(t)|
 \le Cr^2\sup_{0\le t\le r}|\eta(t)|.
\]
For \(Cr^2<1\), we obtain \(\eta=0\) on \([0,r]\).
Uniqueness for \eqref{eq:radial-kernel} away from zero then
gives \(\eta\equiv0\), a contradiction.
Thus \(\eta(0)\ne0\), and we normalize \(\eta(0)=1\).

Since \(\eta\in\D_{\rm rad}\), the Cauchy--Schwarz inequality gives
\[
 \begin{aligned}
 \int_r^\infty|\eta'(t)|\,dt
 &\le
 \left(\int_r^\infty t^{N-1}|\eta'(t)|^2\,dt\right)^{1/2}
 \left(\int_r^\infty t^{1-N}\,dt\right)^{1/2}\\
 &\le Cr^{-(N-2)/2}\|\nabla\eta\|_2
 \qquad(r>0).
 \end{aligned}
\]
The integrability of \(\eta'\) gives a finite limit of \(\eta\)
at infinity, which is zero since \(\eta\in L^{2^*}(\R^N)\). Hence
\[
 \lim_{r\to\infty}\eta(r)=0,
 \qquad
 |\eta(r)|\le\int_r^\infty|\eta'(t)|\,dt
 \le Cr^{-(N-2)/2}\|\nabla\eta\|_2.
\]
Together with local continuity, this bound gives \(\eta\in L^\infty\).

By Lemma~\ref{lem:negative-direction}, \(Z\in\D_{\rm rad}\),
\(L_QZ=-pQ^{p-1}\), and \(\B_Q\) is continuous on \(\D\). Thus
we may test \(L_Q\eta=0\) with \(Z\). The symmetry of
\(\B_Q\) gives
\begin{equation}\label{eq:orth}
 \begin{aligned}
 0=\B_Q(\eta,Z)
 &=\B_Q(Z,\eta)\\
 &=-p|\mathbb S^{N-1}|
   \int_0^\infty r^{N-1}Q(r)^{p-1}\eta(r)\,dr.
 \end{aligned}
\end{equation}
The integral in \eqref{eq:orth} is absolutely convergent because
\(\eta\) is bounded and \(M<\infty\) by \eqref{eq:fast1}.
If \(\eta\ge0\) everywhere, then \(\eta(0)=1\) and \(Q>0\)
would make the integral in \eqref{eq:orth} positive. Thus \(\eta\)
takes negative values and has at least one zero in \((0,\infty)\).

Every zero of \(\eta\) in \((0,\infty)\) is simple: otherwise,
\(\eta(r)=\eta'(r)=0\) at some \(r>0\), and uniqueness for
\eqref{eq:radial-kernel} would give \(\eta\equiv0\).
We now rule out a second zero. Suppose the first two zeros are
\(r_1<r_2\). Keep \(\eta\) on each of the first two intervals
and extend the resulting functions by zero:
\[
 w_1(r)=
 \begin{cases}
  \eta(r),&0\le r\le r_1,\\
  0,&r>r_1,
 \end{cases}
 \qquad
 w_2(r)=
 \begin{cases}
  \eta(r),&r_1\le r\le r_2,\\
  0,&0\le r<r_1\text{ or }r>r_2.
 \end{cases}
\]
Since \(\eta(r_1)=\eta(r_2)=0\), both \(w_1,w_2\) belong to
\(\D_{\rm rad}\), with gradients equal to \(\nabla\eta\)
on their defining intervals and zero elsewhere.
Testing \(L_Q\eta=0\) with \(w_1\) and \(w_2\) gives
\begin{equation}\label{eq:nodal-pieces}
 \B_Q(w_i,w_i)=\B_Q(\eta,w_i)=0\quad(i=1,2),
 \qquad \B_Q(w_1,w_2)=0.
\end{equation}
The cross term is zero because the two functions have disjoint
supports apart from their common boundary.

The linear functional \(w\mapsto\B_Q(w,Z)\) has a nonzero
kernel on the two-dimensional space \(\operatorname{span}\{w_1,w_2\}\).
Choose
\begin{equation}\label{eq3:20261107}
 0\ne w\in\operatorname{span}\{w_1,w_2\},
 \qquad \B_Q(w,Z)=0.
\end{equation}
Equation \eqref{eq:nodal-pieces} gives \(\B_Q(w,w)=0\).

The assumption \(\ind(Q)=1\) and the inequality
\(\B_Q(Z,Z)<0\) imply
\begin{equation}\label{eq2:20261107}
 \B_Q(h,h)\ge0
 \qquad\text{whenever }h\in\D_{\rm rad}
 \text{ and }\B_Q(h,Z)=0.
\end{equation}
Otherwise, \(\B_Q(h,Z)=0\) and \(\B_Q(h,h)<0\) would make
\(\B_Q\) negative definite on \(\operatorname{span}\{h,Z\}\),
contrary to \(\ind(Q)=1\).

Let \(h\in\D_{\rm rad}\) satisfy \(\B_Q(h,Z)=0\).
Then \(\B_Q(w+th,Z)=0\) for every \(t\in\R\).
Applying \eqref{eq2:20261107} to \(w+th\), we obtain
\[
 0\le\B_Q(w+th,w+th)
 =2t\B_Q(w,h)+t^2\B_Q(h,h)
 \qquad(t\in\R).
\]
Since \(\B_Q(w,w)=0\), nonnegativity for both signs of \(t\)
forces the linear coefficient to vanish. Hence
\begin{equation}\label{eq4:20261107}
 \B_Q(w,h)=0
 \qquad\text{whenever }h\in\D_{\rm rad}
 \text{ and }\B_Q(h,Z)=0.
\end{equation}
For an arbitrary \(h\in\D_{\rm rad}\), set
\[
 \widetilde h
 =h-\frac{\B_Q(h,Z)}{\B_Q(Z,Z)}Z.
\]
Then \(\B_Q(\widetilde h,Z)=0\), and \eqref{eq4:20261107}
gives \(\B_Q(w,\widetilde h)=0\). Combining this equality with
\eqref{eq3:20261107}, we obtain
\[
 \B_Q(w,h)
 =\B_Q(w,\widetilde h)
 +\frac{\B_Q(h,Z)}{\B_Q(Z,Z)}\B_Q(w,Z)=0
 \qquad\text{for every }h\in\D_{\rm rad}.
\]
Taking \(h(x)=\varphi(|x|)\) with
\(\varphi\in C_c^\infty((0,\infty))\) in \(\B_Q(w,h)=0\),
we obtain
\[
 \int_0^\infty r^{N-1}(1-q^2)^{-3/2}w'\varphi'\,dr
 =(p-1)\int_0^\infty r^{N-1}Q^{p-2}w\varphi\,dr.
\]
Thus \(w\) satisfies \eqref{eq:radial-kernel} weakly on
\((0,\infty)\). The coefficients are \(C^1\) on compact
subintervals and the leading coefficient is positive, so \(w\)
is a classical solution, including at \(r_1\) and \(r_2\).

Since \(w=0\) for \(r>r_2\), we have \(w(R)=w'(R)=0\)
at any \(R>r_2\). Uniqueness gives \(w\equiv0\), contrary to
\eqref{eq3:20261107}. Hence \(\eta\) has exactly one zero
\(r_1>0\). Since \(\eta(0)=1\) and \(\eta\) takes negative values,
\[
 \eta(r)>0\quad(0\le r<r_1),
 \qquad
 \eta(r)<0\quad(r>r_1).
\]

\medskip\noindent
Step 2. We prove two integral identities and find the sign of \(j\).

Since \(j(0)=0\), we have
\(|j(r)|\le\int_0^r|j'(t)|\,dt\).
Also, \(Q'=-q\) and \(Q(\infty)=0\) give
\(Q(t)=\int_t^\infty q(r)\,dr\).
Tonelli's theorem and \eqref{eq:j} give
\begin{equation}\label{eq:tonelli}
 \begin{aligned}
 \int_0^\infty q(r)\int_0^r|j'(t)|\,dt\,dr
 &=\int_0^\infty |j'(t)|\int_t^\infty q(r)\,dr\,dt\\
 &=\int_0^\infty Q(t)|j'(t)|\,dt\\
 &=(p-1)\int_0^\infty t^{N-1}Q(t)^{p-1}|\eta(t)|\,dt\\
 &\le(p-1)M\|\eta\|_\infty<\infty.
 \end{aligned}
\end{equation}
Here \(M\) is defined in \eqref{eq:fast1}, and \(\eta\) is
bounded by Step 1. In particular,
\[
 \int_0^\infty q(r)|j(r)|\,dr
 \le\int_0^\infty q(r)\int_0^r|j'(t)|\,dt\,dr<\infty.
\]
By \eqref{eq:tonelli}, \(q(r)j'(t)\) is absolutely integrable
on \(0<t<r<\infty\). Changing the order of integration and
using \eqref{eq:j} and \eqref{eq:orth}, we obtain
\begin{equation}\label{eq:momfirst}
 \begin{aligned}
 \int_0^\infty q(r)j(r)\,dr
 &=\int_0^\infty q(r)\int_0^r j'(t)\,dt\,dr\\
 &=\int_0^\infty j'(t)\int_t^\infty q(r)\,dr\,dt\\
 &=\int_0^\infty Q(t)j'(t)\,dt\\
 &=(p-1)\int_0^\infty t^{N-1}Q(t)^{p-1}\eta(t)\,dt=0.
 \end{aligned}
\end{equation}

Next, \eqref{eq:radial-flux-def}, \eqref{eq:flux2}, and
\eqref{eq:j} give
\[
 m'=r^{N-1}Q^{p-1},
 \qquad -m\eta'=q(1-q^2)j.
\]
For integration by parts, we use
\[
 \int_0^\infty|m'\eta|\,dr\le M\|\eta\|_\infty,
 \qquad
 \int_0^\infty|m\eta'|\,dr
 \le\int_0^\infty q|j|\,dr<\infty.
\]
At zero, \(m(0)=0\) and \(\eta\) is bounded. At infinity,
\(m(r)\to M<\infty\) and \(\eta(r)\to0\).
Thus \([m\eta]_0^\infty=0\).
Integrating by parts and using \eqref{eq:orth}, we obtain
\begin{equation}\label{eq:momsecond}
 \begin{aligned}
 0
 &=\int_0^\infty r^{N-1}Q^{p-1}\eta\,dr
 =\int_0^\infty m'\eta\,dr\\
 &=-\int_0^\infty m\eta'\,dr
 =\int_0^\infty q(1-q^2)j\,dr.
 \end{aligned}
\end{equation}
Subtracting \eqref{eq:momsecond} from \eqref{eq:momfirst}
gives the second identity below; the first is \eqref{eq:momfirst}:
\begin{equation}\label{eq:twomom}
 \int_0^\infty qj\,dr=0,
 \qquad
 \int_0^\infty q^3j\,dr=0.
\end{equation}
Both integrals are absolutely convergent since \(0<q<1\) and
\(\int_0^\infty q|j|\,dr<\infty\).

The sign of \(\eta\) from Step 1 and the equation for \(j'\)
in \eqref{eq:j} imply
\[
 j'(r)>0\quad(0<r<r_1),
 \qquad
 j'(r)<0\quad(r>r_1).
\]
Since \(j(0)=0\), we have \(j>0\) on \((0,r_1]\).
If \(j\) never became negative, then \(q>0\) would give
\(\int_0^\infty qj\,dr>0\), contrary to \eqref{eq:twomom}.
Since \(j\) strictly decreases after \(r_1\), the flux \(j\)
has exactly one zero \(r_j>r_1\). Thus
\[
 j(r)>0\quad(0<r<r_j),
 \qquad
 j(r)<0\quad(r>r_j).
\]
It follows that
\begin{equation}\label{eq:partial}
 \int_0^r q(t)j(t)\,dt>0
 \qquad\text{for every }r>0.
\end{equation}
Indeed, for \(0<r\le r_j\), the integrand is positive on \((0,r)\).
For \(r>r_j\), the first identity in \eqref{eq:twomom} gives
\[
 \int_0^r q(t)j(t)\,dt
 =-\int_r^\infty q(t)j(t)\,dt>0.
\]

\medskip\noindent
Step 3. We show that \(q\) first increases and then decreases.

To determine the sign of \(q'\), set
\[
 \psi(r)=(N-1)m(r)-rm'(r).
\]
Differentiating \eqref{eq:flux2} with respect to \(r\) gives
\[
 q'=(1-q^2)^{3/2}
     \left(\frac{m'}{r^{N-1}}-\frac{(N-1)m}{r^N}\right).
\]
Also, since \(m'=r^{N-1}Q^{p-1}\) and \(Q'=-q\),
\[
 \begin{aligned}
 \psi'
 &=(N-2)m'-rm''\\
 &=-r^{N-1}Q^{p-1}
   +(p-1)r^NQ^{p-2}q.
 \end{aligned}
\]
Consequently,
\begin{equation}\label{eq:psidirect}
 q'=-\frac{(1-q^2)^{3/2}}{r^N}\psi,
 \qquad
 \psi'=(p-1)r^{N-1}Q^{p-2}
         \left(rq-\frac{Q}{p-1}\right).
\end{equation}
The first identity in \eqref{eq:psidirect} shows that \(q'\)
and \(\psi\) have opposite signs. To find the sign of \(\psi\),
we compare \(\psi\) with \(j\) using
\[
 W(r)=\eta(r)\psi(r)
      -j(r)\left(rq(r)-\frac{Q(r)}{p-1}\right).
\]
Equations \eqref{eq:j} and \eqref{eq:psidirect} give
\begin{equation}\label{eq:wr-products}
 \begin{aligned}
 \eta'\psi&=jrq',\\
 \eta\psi'&=j'\left(rq-\frac{Q}{p-1}\right),\\
 j'\psi-j\psi'&=(p-1)r^{N-1}Q^{p-2}W.
 \end{aligned}
\end{equation}
Differentiating \(W\) and using \eqref{eq:wr-products}, we obtain
\[
 \begin{aligned}
 W'
 &=\eta'\psi+\eta\psi'
   -j'\left(rq-\frac{Q}{p-1}\right)
   -j\left(rq'+\frac{p}{p-1}q\right)\\
 &=-\frac{p}{p-1}qj.
 \end{aligned}
\]
Near zero, \eqref{eq:flux1} gives \(m=O(r^N)\) and
\(m'=O(r^{N-1})\), while \eqref{eq:flux2} gives \(q=O(r)\).
Also, \eqref{eq:j} gives
\(j=O(r^N)\), and \(\eta\) is bounded. Therefore,
\(\psi=O(r^N)\) and \(W(r)\to0\) as \(r\downarrow0\).
Integrating the equation for \(W'\) and using \eqref{eq:partial},
we obtain
\begin{equation}\label{eq:wr}
 W(r)=-\frac{p}{p-1}\int_0^r q(t)j(t)\,dt<0
 \qquad(r>0).
\end{equation}

\begin{figure}[!htbp]
\centering
\includegraphics[width=.90\linewidth]{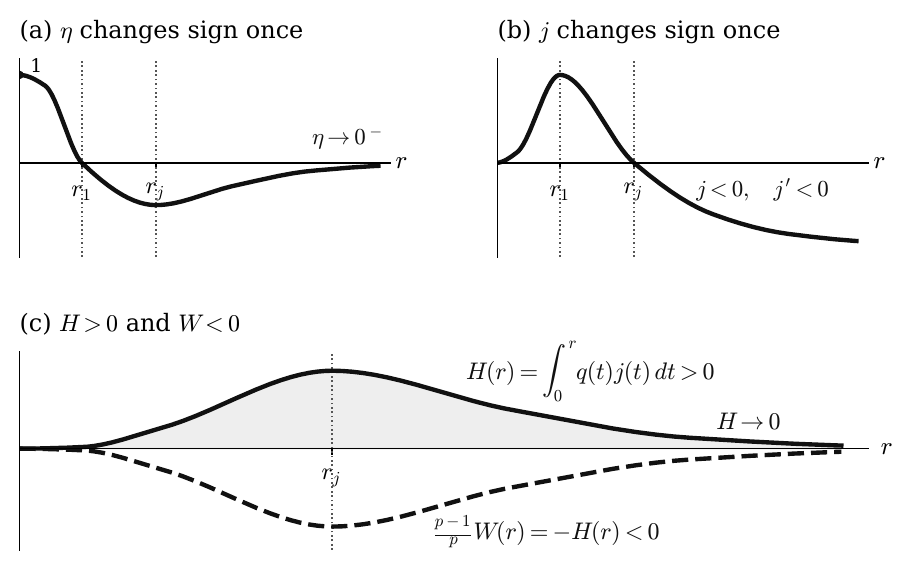}
\caption{The signs of $\eta$, $j$, and $W$.
The function $\eta$ vanishes at $r_1$ and reaches its minimum at
$r_j$, where $j$ vanishes. The partial integral
$H(r)=\int_0^r qj$ stays positive and tends to zero. Hence
$W=-pH/(p-1)<0$. The curves are sketches; their vertical scales
are independent.}
\label{fig:rad-signs}
\end{figure}

Figure~\ref{fig:rad-signs} summarizes the signs of \(\eta,j,W\).

Since \(r_j>r_1\), we have \(\eta(r_j)<0\).
Using \(j(r_j)=0\) in the definition of \(W\), we obtain
\[
 \psi(r_j)=\frac{W(r_j)}{\eta(r_j)}>0.
\]
Since \(\psi\) is continuous and \(j\) changes from positive
to negative at \(r_j\),
\[
 \lim_{r\to r_j^-}\frac{\psi(r)}{j(r)}=+\infty,
 \qquad
 \lim_{r\to r_j^+}\frac{\psi(r)}{j(r)}=-\infty.
\]
At zero, \eqref{eq:flux1} gives
\[
 \lim_{r\downarrow0}\frac{m(r)}{r^N}
 =\frac{Q(0)^{p-1}}{N},
 \qquad
 \lim_{r\downarrow0}\frac{\psi(r)}{r^N}
 =-\frac{Q(0)^{p-1}}{N}<0.
\]
At infinity, \eqref{eq:Q-decay} and \((N-2)(p-1)>N\) imply
\(r^NQ(r)^{p-1}\to0\).
Together with \eqref{eq:fast1}, this yields
\[
 \psi(r)=(N-1)m(r)-r^NQ(r)^{p-1}
 \longrightarrow(N-1)M>0.
\]

For \(r\ne r_j\), equations \eqref{eq:wr-products} and
\eqref{eq:wr} give
\begin{equation}\label{eq:ratio}
 \left(\frac{\psi}{j}\right)'
 =-\frac{j'\psi-j\psi'}{j^2}
 =-\frac{(p-1)r^{N-1}Q^{p-2}}{j^2}W>0.
\end{equation}
On \((0,r_j)\), \(j>0\). The quotient \(\psi/j\) is continuous
and strictly increasing, negative near zero, and positive near
\(r_j\). Hence \(\psi/j\) has exactly one zero
\(r_c\in(0,r_j)\).
Hence
\[
 \psi<0\quad\text{on }(0,r_c),
 \qquad
 \psi>0\quad\text{on }(r_c,r_j].
\]
On \((r_j,\infty)\), the quotient \(\psi/j\) is strictly increasing.
For any \(r>r_j\), choose \(s>r\) with \(\psi(s)>0\), which
is possible since \(\psi(s)\to(N-1)M>0\). Since \(j(s)<0\),
\[
 \frac{\psi(r)}{j(r)}
 <\frac{\psi(s)}{j(s)}<0.
\]
Since \(j(r)<0\), we obtain \(\psi(r)>0\) for every \(r>r_j\).
The equation for \(q'\) in \eqref{eq:psidirect} now gives
\[
 q'(r)>0\quad(0<r<r_c),
 \qquad
 q'(r)<0\quad(r>r_c).
\]

\begin{figure}[!htbp]
\centering
\includegraphics[width=.90\linewidth]{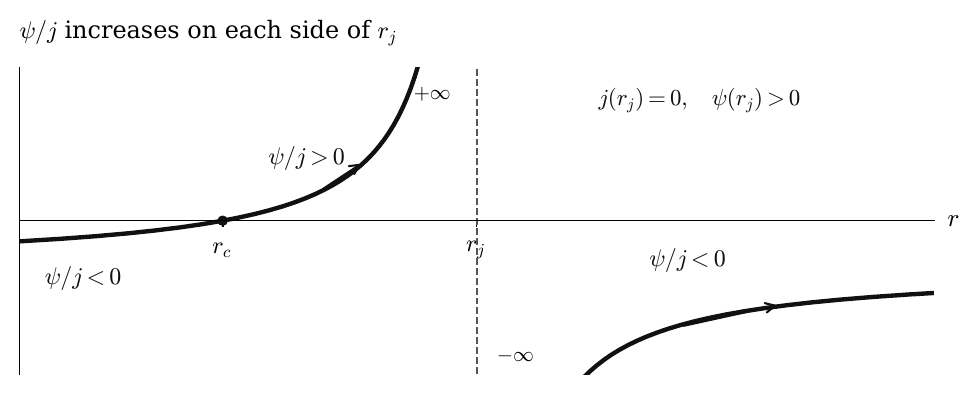}
\caption{The quotient $\psi/j$ increases on each side of $r_j$.
The left part crosses zero once, at $r_c$. The right part stays
negative. Since $j$ changes sign at $r_j$, $\psi$ is positive
for all $r>r_c$. Thus $q$ increases before $r_c$ and decreases
after it. The quotient is undefined at $r_j$, but $\psi(r_j)>0$.}
\label{fig:rad-ratio}
\end{figure}

Figure~\ref{fig:rad-ratio} shows the two parts of \(\psi/j\);
at \(r_j\), we have already proved \(\psi(r_j)>0\).

\medskip\noindent
Step 4. We use the monotonicity of \(q\) to show that an
integral forced to be zero is strictly positive.

Set \(a=q(r_j)\). Since \(r_c<r_j\), we have
\(0<a<q(r_c)<1\). Also, \(q(0)=0\), and \eqref{eq:fast2} gives
\(q(r)\to0\) as \(r\to\infty\).
The monotonicity of \(q\) from Step 3 therefore gives a unique
\(r_a\in(0,r_c)\) such that \(q(r_a)=a\). In particular,
\[
 \{r>0:0<q(r)<a\}
 =(0,r_a)\cup(r_j,\infty),
 \qquad
 \{r>0:q(r)>a\}=(r_a,r_j).
\]
Write \(I_-=(0,r_a)\) and \(I_+=(r_j,\infty)\).
The signs on the two intervals are
\[
 \begin{array}{c|ccc}
  &j&\psi&q'\\ \hline
 I_-&>0&<0&>0\\
 I_+&<0&>0&<0.
 \end{array}
\]

The identities \eqref{eq:twomom} give
\(\int_0^\infty q(q^2-\lambda)j\,dr=0\) for every constant
\(\lambda\). We choose \(0<\lambda<a^2\) so that the integral
is positive on each of \((0,r_a)\), \((r_a,r_j)\), and
\((r_j,\infty)\). On \((r_a,r_j)\), positivity follows from
\(q>a\) and \(j>0\). To handle the other two intervals, choose
\begin{equation}\label{eq:cancel}
 \lambda=1-\sqrt{1-a^2}
 =\frac{a^2}{1+\sqrt{1-a^2}},
 \qquad 0<\lambda<a^2.
\end{equation}
The choice \eqref{eq:cancel} gives the cancellation
\begin{equation}\label{eq:cancellation-integral}
 \begin{aligned}
 \int_0^a\frac{x(x^2-\lambda)}{(1-x^2)^{3/2}}\,dx
 &=\left[\frac{1-\lambda}{\sqrt{1-x^2}}
          +\sqrt{1-x^2}\right]_{x=0}^{x=a}=0.
 \end{aligned}
\end{equation}

On \(I_-\), \(q\) increases from \(0\) to \(a\); on \(I_+\),
\(q\) decreases from \(a\) to \(0\). Thus there are unique
\(R_-\in I_-\) and \(R_+\in I_+\) with
\(q(R_-)=q(R_+)=\sqrt\lambda\). We write \(R\) for the point
in the interval under consideration.

\begin{figure}[!htbp]
\centering
\includegraphics[width=.90\linewidth]{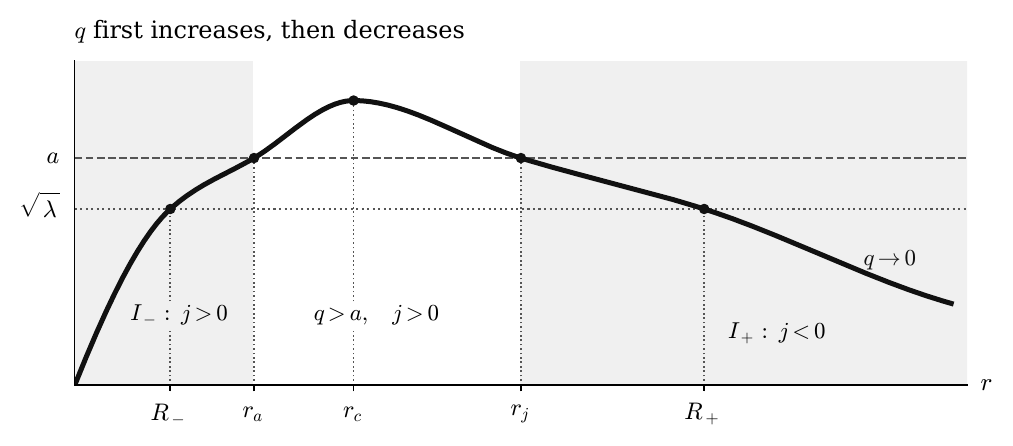}
\caption{The levels $a$ and $\sqrt\lambda$ on the graph of $q$.
We have $q(r_a)=q(r_j)=a$ and
$q(R_-)=q(R_+)=\sqrt\lambda$. Thus
$0<R_-<r_a<r_c<r_j<R_+$. On the middle interval,
$q>a$ and $j>0$. The curve shows only these signs and positions.}
\label{fig:rad-branches}
\end{figure}

Figure~\ref{fig:rad-branches} shows the levels \(a\) and \(\sqrt\lambda\).

On each outer interval, we introduce a factor \(F\) so that
we can use \eqref{eq:cancellation-integral} in
\(\int_I q(q^2-\lambda)j\,dr\). For
\(I\in\{I_-,I_+\}\), set
\[
 F(r)=-\frac{r^N j(r)}{\psi(r)}.
\]
The signs of \(j\) and \(\psi\) give \(F>0\).
Using \eqref{eq:wr-products} and \eqref{eq:wr}, we obtain
\begin{equation}\label{eq:monotonefactor}
 \begin{aligned}
 F'(r)
 &=Nr^{N-1}\frac{-j}{\psi}
   -r^N\frac{j'\psi-j\psi'}{\psi^2}\\
 &=Nr^{N-1}\frac{-j}{\psi}
   -\frac{(p-1)r^{2N-1}Q^{p-2}}{\psi^2}W>0.
 \end{aligned}
\end{equation}
Both terms are positive because \(-j/\psi>0\) and \(W<0\).
Thus \(F\) is strictly increasing. The equation for \(q'\)
in \eqref{eq:psidirect} also gives
\[
 F(r)=\frac{(1-q(r)^2)^{3/2}j(r)}{q'(r)}
 \qquad(r\in I).
\]

On either outer interval \(I\), the substitution \(x=q(r)\)
and \eqref{eq:cancellation-integral} give
\begin{equation}\label{eq:cancel-I}
 \begin{aligned}
 \int_I\frac{q(q^2-\lambda)q'}{(1-q^2)^{3/2}}\,dr
 &=
 \begin{cases}
 \displaystyle\int_0^a
   \frac{x(x^2-\lambda)}{(1-x^2)^{3/2}}\,dx,&I=I_-,\\[6pt]
 \displaystyle\int_a^0
   \frac{x(x^2-\lambda)}{(1-x^2)^{3/2}}\,dx,&I=I_+,
 \end{cases}\\
 &=0.
 \end{aligned}
\end{equation}
The integrals are absolutely convergent, because
\[
 \int_I\frac{q|q^2-\lambda|\,|q'|}{(1-q^2)^{3/2}}\,dr
 =\int_0^a\frac{x|x^2-\lambda|}{(1-x^2)^{3/2}}\,dx<\infty.
\]
Also, \(0<q<1\) and \(\int_0^\infty q|j|\,dr<\infty\) imply
\[
 \int_I|q(q^2-\lambda)j|\,dr
 \le(1+\lambda)\int_I q|j|\,dr<\infty.
\]
We may therefore subtract \(F(R)\) times \eqref{eq:cancel-I}.
Using the definition of \(F\), we obtain
\begin{equation}\label{eq:onesign}
 \begin{aligned}
 \int_I q(q^2-\lambda)j\,dr
 &=\int_I F(r)
   \frac{q(q^2-\lambda)q'}{(1-q^2)^{3/2}}\,dr\\
 &=\int_I\bigl(F(r)-F(R)\bigr)
   \frac{q(q^2-\lambda)q'}{(1-q^2)^{3/2}}\,dr>0.
 \end{aligned}
\end{equation}
We now check the strict inequality in \eqref{eq:onesign}. Since \(F\) is
strictly increasing, \(F(r)-F(R)\) is negative before \(R\)
and positive after \(R\). On \(I_-\), \(q\) increases and
\(q'>0\), so \((q^2-\lambda)q'\) has the same signs.
On \(I_+\), \(q^2-\lambda\) is positive before \(R\) and
negative after \(R\), but \(q'<0\). The product
\((q^2-\lambda)q'\) therefore has the same signs as on \(I_-\).
The factor \(q(1-q^2)^{-3/2}\) is positive. Thus the final
integrand in \eqref{eq:onesign} is positive at every point
of \(I\) except \(R\).

\begin{figure}[!htbp]
\centering
\includegraphics[width=.90\linewidth]{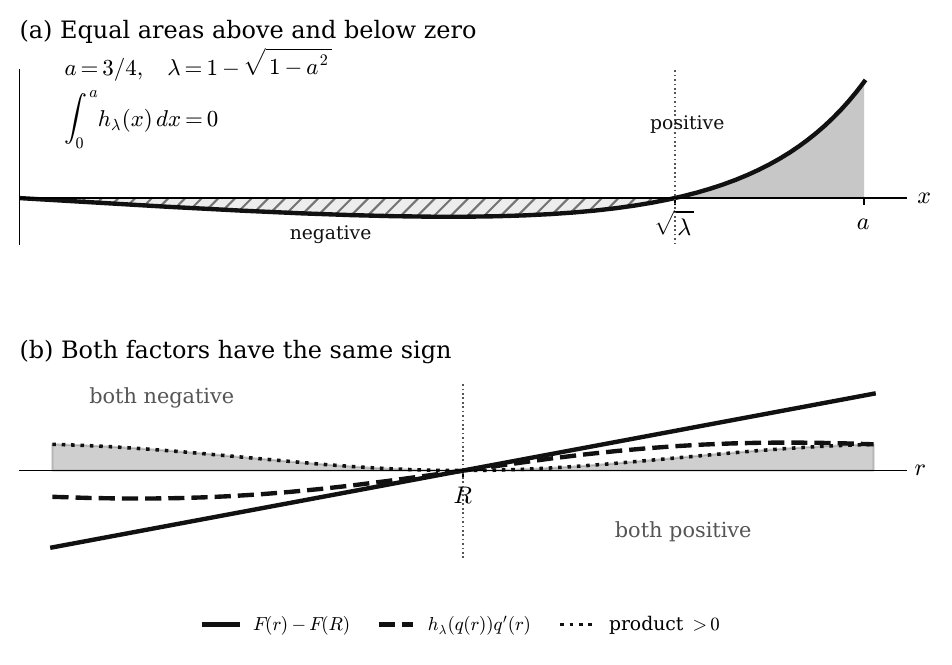}
\caption{Why the final integral is positive.
In (a), $h_\lambda(x)=x(x^2-\lambda)/(1-x^2)^{3/2}$ is plotted
with $a=3/4$ and $\lambda=1-\sqrt{1-a^2}$; the two signed areas
cancel. In (b), both $F(r)-F(R)$ and $h_\lambda(q(r))q'(r)$
are negative before $R$ and positive after it, so their product
is positive except at $R$. Panel (b) shows signs only.}
\label{fig:rad-weight}
\end{figure}

Figure~\ref{fig:rad-weight} shows the matching signs of
\(F(r)-F(R)\) and \((q^2-\lambda)q'\).

On the middle interval \((r_a,r_j)\), we have
\(q>a>\sqrt\lambda\) and \(j>0\), so
\[
 q(q^2-\lambda)j>0.
\]
Adding the three positive integrals and using \eqref{eq:twomom},
we obtain
\[
 0<\int_0^\infty q(q^2-\lambda)j\,dr
 =\int_0^\infty q^3j\,dr
  -\lambda\int_0^\infty qj\,dr=0.
\]
This contradiction proves \(\eta=0\).
\end{proof}
\FloatBarrier

\subsection{Spherical harmonics}

\begin{proof}[Proof of Theorem~\ref{th:nondegeneracy}]
By Lemma~\ref{lem:radial-symmetry}, translate \(Q\) so that
\(Q\) is radial about the origin. Lemma~\ref{lem:ground-morse}
gives \(\ind(Q)=1\). Lemma~\ref{lem:rad} therefore gives
\[
 \ker L_Q\cap\D_{\rm rad}=\{0\}.
\]

Let \(\phi\in\ker_{\D(\R^N)}L_Q\).
Write \(x=r\omega\), where \(r=|x|>0\) and
\(\omega\in\mathbb S^{N-1}\).
We use \(\nabla_{\mathbb S}\) and \(\Delta_{\mathbb S}\)
for the gradient and Laplacian on the unit sphere, acting
on the variable \(\omega\).

For each \(l\ge0\), let \(d_l\) be the dimension of the space
of real-valued spherical harmonics of degree \(l\), and choose
an orthonormal basis \(\{Y_{l,k}:1\le k\le d_l\}\).
The basis functions satisfy
\begin{equation}\label{eq:spherical-eigenvalues}
 -\Delta_{\mathbb S}Y_{l,k}=\lambda_lY_{l,k},
 \qquad
 \lambda_l=l(l+N-2).
\end{equation}
Define
\begin{equation}\label{eq:spherical-coefficients}
 h_{l,k}(r)
 =\int_{\mathbb S^{N-1}}
   \phi(r\omega)Y_{l,k}(\omega)\,d\omega.
\end{equation}
The spherical harmonics form a complete orthonormal basis of
\(L^2(\mathbb S^{N-1})\). Hence
\begin{equation}\label{eq:spherical-expansion}
 \phi(r\omega)
 =\sum_{l=0}^{\infty}\sum_{k=1}^{d_l}
   h_{l,k}(r)Y_{l,k}(\omega),
\end{equation}
with convergence in \(L^2(\mathbb S^{N-1})\) for almost every
\(r>0\).

Fix \(l,k\), and write \(Y_l=Y_{l,k}\) and \(h=h_{l,k}\).
By Cauchy--Schwarz and \(\int_{\mathbb S^{N-1}}Y_l^2\,d\omega=1\),
\[
 |h'(r)|^2
 \le\int_{\mathbb S^{N-1}}
       |\partial_r\phi(r\omega)|^2\,d\omega.
\]
For \(l\ge1\), integration by parts on the sphere and
\eqref{eq:spherical-eigenvalues} give
\[
 \lambda_l h
 =\int_{\mathbb S^{N-1}}
   \nabla_{\mathbb S}\phi\cdot\nabla_{\mathbb S}Y_l\,d\omega,
 \qquad
 \int_{\mathbb S^{N-1}}|\nabla_{\mathbb S}Y_l|^2\,d\omega
 =\lambda_l.
\]
Another application of Cauchy--Schwarz gives
\[
 \lambda_l h^2
 \le\int_{\mathbb S^{N-1}}
       |\nabla_{\mathbb S}\phi|^2\,d\omega.
\]
For \(l=0\), the angular bound is immediate since \(\lambda_0=0\).
Integrating the two bounds gives
\begin{equation}\label{eq:mode-energy}
 \int_0^\infty
 \bigl(r^{N-1}|h'|^2+\lambda_l r^{N-3}h^2\bigr)\dr
 \le\|\nabla\phi\|_2^2.
\end{equation}
These estimates extend from smooth functions to
\(\phi\in\D(\R^N)\) by density.

Recall that \(q=-Q'\) and \(\nabla Q=-q(r)\omega\).
The gradient of \(\phi\) satisfies
\[
 \nabla\phi
 =\partial_r\phi\,\omega+\frac1r\nabla_{\mathbb S}\phi,
 \qquad
 \nabla_{\mathbb S}\phi\cdot\omega=0.
\]
Substituting into the gradient terms of \(L_Q\), we obtain
\[
 \begin{aligned}
 \frac{\nabla\phi}{\sqrt{1-|\nabla Q|^2}}
 +\frac{(\nabla Q\cdot\nabla\phi)\nabla Q}
       {(1-|\nabla Q|^2)^{3/2}}
=
 \frac{\partial_r\phi}{(1-q^2)^{3/2}}\,\omega
 +\frac{\nabla_{\mathbb S}\phi}{r\sqrt{1-q^2}}.
 \end{aligned}
\]
Thus \(L_Q\phi=0\) takes the form
\begin{equation}\label{eq:polar-linear}
 -\frac1{r^{N-1}}\partial_r
 \left(r^{N-1}(1-q^2)^{-3/2}\partial_r\phi\right)
 -\frac{(1-q^2)^{-1/2}}{r^2}\Delta_{\mathbb S}\phi
 =(p-1)Q^{p-2}\phi.
\end{equation}

For \(\xi\in C_c^\infty(0,\infty)\), test
\eqref{eq:polar-linear} with \(\xi(r)Y_l(\omega)\).
By \eqref{eq:spherical-coefficients} and integration by parts
on the sphere,
\[
 \int_{\mathbb S^{N-1}}\partial_r\phi\,Y_l\,d\omega=h'(r),
 \qquad
 \int_{\mathbb S^{N-1}}
 \nabla_{\mathbb S}\phi\cdot\nabla_{\mathbb S}Y_l\,d\omega
 =\lambda_l h(r).
\]
Using \(dx=r^{N-1}\,dr\,d\omega\), we obtain
\begin{equation}\label{eq:mode-weak}
 \begin{aligned}
 \int_0^\infty \Big(
 &r^{N-1}(1-q^2)^{-3/2}h'\xi'\\
 &+\big[\lambda_l r^{N-3}(1-q^2)^{-1/2}
       -(p-1)r^{N-1}Q^{p-2}\big]h\xi
 \Big)\,dr=0.
 \end{aligned}
\end{equation}

We first consider the degree-zero part of \(\phi\), which is
its spherical average:
\[
 \phi_0(x)
 =\frac1{|\mathbb S^{N-1}|}
   \int_{\mathbb S^{N-1}}\phi(|x|\omega)\,d\omega.
\]
H\"older's inequality and \eqref{eq:mode-energy} give
\(\phi_0\in\D_{\rm rad}\).
Since \(Q\) is radial, averaging the weak equation
\(L_Q\phi=0\) over rotations gives \(L_Q\phi_0=0\) on \(\R^N\).
Lemma~\ref{lem:rad} gives \(\phi_0=0\), so \(h_{0,1}=0\).

We now consider \(l\ge1\). Define
\begin{equation}\label{eq:mode}
 \begin{aligned}
 B_l(v)=\int_0^\infty\Big(
 &r^{N-1}(1-q^2)^{-3/2}|v'|^2\\
 &+\lambda_l r^{N-3}(1-q^2)^{-1/2}v^2
 -(p-1)r^{N-1}Q^{p-2}v^2
 \Big)\dr.
 \end{aligned}
\end{equation}
Differentiating \eqref{eq:pde} with respect to \(x_i\) gives
\begin{equation}\label{eq:translation-kernel}
 L_Q(\partial_{x_i}Q)=0,
 \qquad
 \partial_{x_i}Q(x)=-q(r)\frac{x_i}{r}.
\end{equation}
Regularity at zero and the bounds in \eqref{eq:fast2}
give \(\partial_{x_i}Q\in\D(\R^N)\).
Since \(x_i/r\) is a spherical harmonic of degree one,
\(q\) satisfies \eqref{eq:mode-weak} with \(l=1\), for which
\(\lambda_1=N-1\).

Let \(v\in H^1(0,\infty)\) have compact support in
\((0,\infty)\). Since \(q>0\), we may use \(\xi=v^2/q\)
in \eqref{eq:mode-weak} with \(h=q\) and \(l=1\), by density
on an interval containing the support of \(v\). We obtain
\[
 \begin{aligned}
 &\int_0^\infty
 r^{N-1}(1-q^2)^{-3/2}q'(v^2/q)'\dr\\
 &\qquad=
 \int_0^\infty
 \bigl((p-1)r^{N-1}Q^{p-2}
       -\lambda_1 r^{N-3}(1-q^2)^{-1/2}\bigr)v^2\dr.
 \end{aligned}
\]
Combining this equality with
\[
 |v'|^2-q'(v^2/q)'=q^2|(v/q)'|^2
\]
and the definition \eqref{eq:mode}, we obtain
\begin{equation}\label{eq:modeidentity}
 \begin{aligned}
 B_l(v)
 &=\int_0^\infty
 r^{N-1}(1-q^2)^{-3/2}q^2|(v/q)'|^2\dr\\
 &\quad+(\lambda_l-\lambda_1)
 \int_0^\infty r^{N-3}(1-q^2)^{-1/2}v^2\dr.
 \end{aligned}
\end{equation}

Choose \(0<2\delta<R\) and a smooth cutoff \(\chi\) such that
\[
 0\le\chi\le1,\qquad
 \chi=1\ \text{on }[2\delta,R],\qquad
 \chi=0\ \text{outside }[\delta,2R],\qquad
 |\chi'|\le C/r.
\]
We may use \(\xi=\chi^2h\) in \eqref{eq:mode-weak} by density
on \([\delta,2R]\). Hence
\begin{equation}\label{eq:angular-cut}
 \begin{aligned}
 B_l(\chi h)
 &=\int_0^\infty r^{N-1}(1-q^2)^{-3/2}
   \bigl(|(\chi h)'|^2-h'(\chi^2h)'\bigr)\dr\\
 &=\int_0^\infty r^{N-1}(1-q^2)^{-3/2}
   |\chi'|^2h^2\dr\\
 &\le C\left(\int_\delta^{2\delta}+\int_R^{2R}\right)
   r^{N-3}h^2\dr.
 \end{aligned}
\end{equation}
The constant \(C\) is independent of \(\delta,R\), since
\(\|q\|_\infty<1\).
For \(l\ge1\), \eqref{eq:mode-energy} gives
\[
 \int_0^\infty r^{N-3}h^2\dr<\infty.
\]
Therefore the right-hand side of \eqref{eq:angular-cut}
tends to zero as \(\delta\downarrow0\) and \(R\to\infty\).

Fix a compact interval \(J\subset(0,\infty)\).
For sufficiently small \(\delta\) and sufficiently large \(R\),
we have \(\chi=1\) on \(J\).
Since \(\lambda_l\ge\lambda_1\), \eqref{eq:modeidentity}
and \eqref{eq:angular-cut} give
\begin{equation}\label{eq:mode-local}
 \begin{aligned}
 0\le
 &\int_J r^{N-1}(1-q^2)^{-3/2}q^2|(h/q)'|^2\dr\\
 &+(\lambda_l-\lambda_1)
   \int_J r^{N-3}(1-q^2)^{-1/2}h^2\dr\\
 &\le B_l(\chi h)\longrightarrow0.
 \end{aligned}
\end{equation}
For \(l=1\), \eqref{eq:mode-local} gives \((h/q)'=0\)
on every compact interval in \((0,\infty)\).
Thus
\[
 h_{1,k}(r)=a_kq(r)
\]
for some constant \(a_k\).
For \(l\ge2\), we have \(\lambda_l>\lambda_1\).
The second integral in \eqref{eq:mode-local} then gives
\[
 h_{l,k}=0.
\]

We have also proved \(h_{0,1}=0\).
Thus only the degree-one terms remain in
\eqref{eq:spherical-expansion}.
The degree-one spherical harmonics are spanned by the
coordinate functions
\[
 \omega\longmapsto\omega_i=\frac{x_i}{|x|},
 \qquad i=1,\ldots,N,
 \qquad \omega=\frac{x}{|x|}\in\mathbb S^{N-1}.
\]
Thus \(d_1=N\), and for some constants \(c_1,\ldots,c_N\),
\[
 \begin{aligned}
 \phi(r\omega)
 &=q(r)\sum_{k=1}^{N}a_kY_{1,k}(\omega)\\
 &=q(r)\sum_{i=1}^{N}c_i\omega_i
 =-\sum_{i=1}^{N}c_i\,\partial_{x_i}Q(r\omega).
 \end{aligned}
\]
Every kernel function is therefore a linear combination of
\(\partial_{x_1}Q,\ldots,\partial_{x_N}Q\).
Conversely, \eqref{eq:translation-kernel} and the bounds
\eqref{eq:fast2} show that every translation derivative belongs
to \(\ker_{\D(\R^N)}L_Q\). This proves \eqref{eq:kernel}.
\end{proof}

\section{Continuation in the exponent}\label{sec:cont}
We now show that the number of positive radial finite-energy
solutions with radial Morse index one does not change with \(p\).

\subsection{Compactness of radial solutions}\label{sec:compact}

We first prove uniform bounds for radial solutions. The bounds
will be used in the compactness, continuation, and uniqueness proofs.

\begin{lemma}\label{lem:radial-bounds}
Let \(N\ge3\) and \(p_0>2^*\). For a positive radial finite-energy solution
\(Q\) of \eqref{eq:pde} with exponent \(p\ge p_0\), set
\[
 \alpha=Q(0),\qquad q=-Q',\qquad
 m(r)=\int_0^r t^{N-1}Q(t)^{p-1}\,dt,\qquad
 M=\lim_{r\to\infty}m(r).
\]
There are \(c,C>0\) and \(R\ge1\), depending only on
\(N,p_0\), such that every such solution satisfies
\[
 c\le\alpha\le C,\qquad
 \alpha^p\le Cp^N,\qquad \alpha^{p-2^*}\ge c,
\]
\[
 \|\nabla Q\|_2^2+\|Q\|_p^p+M\le C,
 \qquad
 Q(r)\le Cr^{2-N},\quad q(r)\le Cr^{1-N}\quad(r\ge1).
\]
In addition,
\[
 m(R)\ge\frac34M,\qquad q(r)\le\frac12\quad(r\ge R).
\]

\end{lemma}

\begin{proof}
By \eqref{eq:flux2},
\begin{equation}\label{eq:compact-flux}
 m'(r)=r^{N-1}Q(r)^{p-1},\qquad
 q(r)=\frac{m(r)}{\sqrt{r^{2N-2}+m(r)^2}}.
\end{equation}
In particular, \(m\) is increasing and \(0<q(r)<1\) for
\(r>0\). Since \(Q(r)\to0\) as \(r\to\infty\), we also have
\begin{equation}\label{eq:compact-tail}
 Q(r)=\int_r^\infty q(t)\,dt.
\end{equation}

We first bound \(\alpha\). Since \(|Q'|\le1\),
\(Q(t)\ge\alpha-t\ge\alpha(1-1/p)\) for
\(0\le t\le\alpha/p\). Hence
\begin{equation}\label{eq:core-height}
 \begin{aligned}
 m(\alpha/p)
 &=\int_0^{\alpha/p}t^{N-1}Q(t)^{p-1}\,dt\\
 &\ge [\alpha(1-1/p)]^{p-1}
       \int_0^{\alpha/p}t^{N-1}\,dt\\
 &=\frac{\alpha^{p+N-1}}{Np^N}(1-1/p)^{p-1}
 \ge\frac{\alpha^{p+N-1}}{eNp^N}.
 \end{aligned}
\end{equation}
The last inequality uses \((1-1/p)^{p-1}\ge e^{-1}\).

To bound \(Q(\alpha/p)\) from below in terms of \(m(\alpha/p)\), set
\[
 c_N:=\int_0^\infty\frac{dv}{\sqrt{1+v^{2N-2}}}.
\]
Since the integrand is bounded near zero and is \(O(v^{1-N})\)
at infinity, \(0<c_N<\infty\) for \(N\ge3\). The substitution
\(t=s^{1/(N-1)}v\) gives
\begin{equation}\label{eq:gint}
 \int_0^\infty\frac{s}{\sqrt{t^{2N-2}+s^2}}\,dt
 =c_Ns^{1/(N-1)}\qquad(s>0).
\end{equation}
For fixed \(t>0\), the function
\(s\mapsto s/\sqrt{t^{2N-2}+s^2}\) is increasing.
Since \(m(t)\ge m(\alpha/p)\) for \(t\ge\alpha/p\),
\eqref{eq:compact-flux} and \eqref{eq:compact-tail} give
\[
 \begin{aligned}
 \alpha\ge Q(\alpha/p)
 &\ge\int_{\alpha/p}^\infty
 \frac{m(\alpha/p)}{\sqrt{t^{2N-2}+m(\alpha/p)^2}}\,dt\\
 &\ge c_Nm(\alpha/p)^{1/(N-1)}-\frac{\alpha}{p}.
 \end{aligned}
\]
Here we used \eqref{eq:gint} with \(s=m(\alpha/p)\) and
the bound \(s/\sqrt{t^{2N-2}+s^2}\le1\) on \((0,\alpha/p)\).
Consequently,
\[
 \frac{\alpha^{p+N-1}}{eNp^N}
 \le m(\alpha/p)
 \le \left(\frac{\alpha(1+1/p)}{c_N}\right)^{N-1}.
\]
Hence
\begin{equation}\label{eq:heightupper}
 \alpha^p\le C_0p^N,\qquad \alpha\le C_1,
\end{equation}
where \(C_0>0\) depends only on \(N\), and \(C_1>0\)
depends only on \(N,p_0\). The second bound follows since
\((C_0p^N)^{1/p}\) is bounded for \(p\ge p_0\).

For \(t\ge r>0\), \eqref{eq:compact-flux} and
\(m(t)\ge m(r)\) give
\[
 \begin{aligned}
 q(t)
 &\ge\frac{m(r)}{\sqrt{t^{2N-2}+m(r)^2}}\\
 &=q(r)\frac{\sqrt{r^{2N-2}+m(r)^2}}
                  {\sqrt{t^{2N-2}+m(r)^2}}
 \ge q(r)(r/t)^{N-1}.
 \end{aligned}
\]
Integrating this inequality and using \eqref{eq:compact-tail},
we obtain
\[
 Q(r)\ge q(r)\int_r^\infty(r/t)^{N-1}\,dt
 =\frac{r}{N-2}q(r).
\]
Therefore,
\begin{equation}\label{eq:rq}
 rq(r)\le(N-2)Q(r)\le(N-2)\alpha\le(N-2)C_1.
\end{equation}

We next bound the Dirichlet integral
\begin{equation}\label{eq:compact-dirichlet}
 \int_{\mathbb R^N}|\nabla Q|^2\,dx
 =\|\nabla Q\|_2^2
 =|\mathbb S^{N-1}|\int_0^\infty r^{N-1}q(r)^2\,dr.
\end{equation}
The estimate \(q(r)\le (N-2)C_1/r\) in \eqref{eq:rq} does
not bound \eqref{eq:compact-dirichlet}, since
\(\int_1^\infty r^{N-3}\,dr=\infty\).
We instead bound \(r^{N-1}q(r)^2\) by \(-\mathcal P'(r)\)
for large \(r\), where
\begin{equation}\label{eq:compact-P}
 \mathcal P(r):=
 r^N\left(\frac1{\sqrt{1-q(r)^2}}-1+\frac{Q(r)^p}{p}\right)
 -\frac Np m(r)Q(r).
\end{equation}
The formula for \(q\) in \eqref{eq:compact-flux} implies
\[
 \left((1-q^2)^{-1/2}\right)'
 =q\left(\frac{m}{r^{N-1}}\right)',
 \qquad
 r^{N-1}\big((1-q^2)^{-1/2}-1\big)
 =\frac{qm}{1+\sqrt{1-q^2}}.
\]
The choice of \(\mathcal P\) makes the terms containing
\(Q^p\) and \(qQ^{p-1}\) cancel after differentiation.
Using \(m'=r^{N-1}Q^{p-1}\) and \(Q'=-q\), we obtain
\begin{equation}\label{eq:localP}
 \begin{aligned}
 \mathcal P'(r)
 &=Nr^{N-1}\big((1-q(r)^2)^{-1/2}-1\big)
   -\left(N-1-\frac Np\right)q(r)m(r)\\
 &=q(r)m(r)\left(
 1+\frac Np-
 \frac{N\sqrt{1-q(r)^2}}{1+\sqrt{1-q(r)^2}}
 \right).
 \end{aligned}
\end{equation}
Since \(p_0>2^*\),
\[
 1+\frac N{p_0}-\frac N2<0.
\]
By continuity, there are \(\delta\in(0,1/2]\) and \(c>0\),
depending only on \(N,p_0\), such that
\[
 1+\frac Np-\frac{N\sqrt{1-s^2}}{1+\sqrt{1-s^2}}
 \le-c
 \qquad(p\ge p_0,\ 0\le s\le\delta).
\]
Choose \(R\ge1\), depending only on \(N,p_0\), so that
\((N-2)C_1/R\le\delta\).
By \eqref{eq:rq}, \(q(r)\le\delta\le1/2\) for \(r\ge R\).
Moreover, \eqref{eq:compact-flux} gives
\[
 q(r)m(r)
 =\frac{r^{N-1}q(r)^2}{\sqrt{1-q(r)^2}}
 \ge r^{N-1}q(r)^2.
\]
Integrating \eqref{eq:localP} on \([R,S]\), we obtain
\begin{equation}\label{eq:exterior-dirichlet}
 c\int_R^S r^{N-1}q(r)^2\,dr
 \le\mathcal P(R)-\mathcal P(S)
 \qquad(S>R).
\end{equation}

We estimate \(\mathcal P(R)\) using \(q(R)\le1/2\).
Since \(Q\) is decreasing, \eqref{eq:compact-flux} gives
\[
 \frac{R^N}{N}Q(R)^{p-1}
 \le\int_0^R t^{N-1}Q(t)^{p-1}\,dt
 =m(R)
 =\frac{R^{N-1}q(R)}{\sqrt{1-q(R)^2}}
 \le\frac{R^{N-1}}{\sqrt3}.
\]
Together with \(Q(R)\le C_1\), this gives
\[
 R^NQ(R)^p\le Nm(R)Q(R)
 \le\frac{NC_1R^{N-1}}{\sqrt3}.
\]
Equation \eqref{eq:compact-P} then gives \(|\mathcal P(R)|\le C_2\),
where \(C_2>0\) depends only on \(N,p_0\).

For each fixed solution \(Q\), \eqref{eq:Q-decay},
\eqref{eq:fast1}, and \eqref{eq:fast2} give, as \(r\to\infty\),
\[
 \begin{aligned}
 r^N\big((1-q^2)^{-1/2}-1\big)&=O(r^{2-N}),\\
 r^NQ^p&=O\big(r^{N-(N-2)p}\big),\\
 mQ&=O(r^{2-N}).
 \end{aligned}
\]
The constants here may depend on \(Q,p\). All three terms
tend to zero since \(N\ge3\) and \(p>2^*\), so
\(\mathcal P(S)\to0\). Letting \(S\to\infty\) in
\eqref{eq:exterior-dirichlet} and using \(q\le1\) on \([0,R]\),
we obtain
\[
 \int_0^\infty r^{N-1}q(r)^2\,dr
 \le\frac{R^N}{N}+\frac{C_2}{c}.
\]
Hence \(\|\nabla Q\|_2^2\le C_3\), where \(C_3>0\)
depends only on \(N,p_0\).

We next bound \(\|Q\|_p^p\) and keep \(\alpha\) away from
zero. The identities \eqref{eq:testq} and \eqref{eq:pohoz}
in Lemma~\ref{lem:regularity-energy} give
\[
 \begin{aligned}
 \|\nabla Q\|_2^2
 &\le\|Q\|_p^p\\
 &=\frac{Np}{N+p}
 \int_{\mathbb R^N}\bigl(1-\sqrt{1-|\nabla Q|^2}\bigr)\,dx
 \le N\|\nabla Q\|_2^2.
 \end{aligned}
\]
Here we used \(1-\sqrt{1-s^2}\le s^2\) for \(0\le s\le1\).
Let \(S_N\) be a constant in the Sobolev inequality
\(\|Q\|_{2^*}\le S_N\|\nabla Q\|_2\); \(S_N\) depends
only on \(N\). Since \(Q\le\alpha\),
\[
 \|\nabla Q\|_2^2
 \le\|Q\|_p^p
 \le\alpha^{p-2^*}\|Q\|_{2^*}^{2^*}
 \le S_N^{2^*}\alpha^{p-2^*}\|\nabla Q\|_2^{2^*}.
\]
Since \(0<\|\nabla Q\|_2^2\le C_3\), we obtain
\[
 \alpha^{p-2^*}
 \ge\frac1{S_N^{2^*}C_3^{(2^*-2)/2}}=:c_0>0.
\]
The constant \(c_0\) depends only on \(N,p_0\).
If \(\alpha\le1\), then \(\alpha^{p_0-2^*}\ge c_0\).
Thus \(\alpha\ge\min\{1,c_0^{1/(p_0-2^*)}\}>0\), and
\begin{equation}\label{eq:uniform}
 \|\nabla Q\|_2^2+\|Q\|_p^p\le C,
 \qquad
 \min\{\alpha^{p-2^*},\alpha\}\ge c>0,
\end{equation}
where \(C,c>0\) depend only on \(N,p_0\).

We now prove decay bounds with constants independent of \(Q,p\).
By \eqref{eq:compact-tail}, Cauchy--Schwarz, and
\eqref{eq:uniform},
\begin{equation}\label{eq:slow-decay}
 \begin{aligned}
 Q(r)
 &=\int_r^\infty q(t)\,dt\\
 &\le\left(\int_r^\infty t^{N-1}q(t)^2\,dt\right)^{1/2}
       \left(\int_r^\infty t^{1-N}\,dt\right)^{1/2}
 \le C_4r^{-(N-2)/2},
 \end{aligned}
\end{equation}
where \(C_4>0\) depends only on \(N,p_0\).
For each fixed solution, \(M<\infty\) by \eqref{eq:fast1}.
Since \(m(r)\le M\), \eqref{eq:compact-flux} gives
\(q(r)\le Mr^{1-N}\). Thus
\begin{equation}\label{eq:M-decay}
 Q(r)\le\int_r^\infty Mt^{1-N}\,dt
 =\frac{M}{N-2}r^{2-N}.
\end{equation}

It remains to bound \(M\) uniformly.
Increase \(R\), still depending only on \(N,p_0\), so that
\(C_4t^{-(N-2)/2}\le1\) for \(t\ge R\).
By \eqref{eq:slow-decay}, \(Q(t)\le1\) on \([R,\infty)\).
Since \(p\ge p_0\),
\[
 Q(t)^{p-2}\le Q(t)^{p_0-2}
 \le C_4^{p_0-2}t^{-(N-2)(p_0-2)/2}
 \qquad(t\ge R).
\]
Using the bound for \(Q^{p-2}\) and \eqref{eq:M-decay}, we obtain
\begin{equation}\label{eq:absorbM}
 \begin{aligned}
 M-m(R)
 &=\int_R^\infty t^{N-1}Q(t)^{p-1}\,dt\\
 &\le\frac{M}{N-2}\int_R^\infty tQ(t)^{p-2}\,dt\\
 &\le\frac{C_4^{p_0-2}M}{N-2}
       \int_R^\infty t^{1-(N-2)(p_0-2)/2}\,dt.
 \end{aligned}
\end{equation}
Since \((N-2)(p_0-2)/2>2\), the integral in
\eqref{eq:absorbM} tends to zero as \(R\to\infty\).
Increase \(R\) again, depending only on \(N,p_0\), so that
\[
 \frac{C_4^{p_0-2}}{N-2}
 \int_R^\infty t^{1-(N-2)(p_0-2)/2}\,dt\le\frac14.
\]
The bound \(q(R)\le1/2\) still holds by \eqref{eq:rq}.
Therefore, \eqref{eq:absorbM} gives
\[
 \frac34M\le m(R)
 =\frac{R^{N-1}q(R)}{\sqrt{1-q(R)^2}}
 \le\frac{R^{N-1}}{\sqrt3}.
\]
Combining \(M\le4R^{N-1}/(3\sqrt3)\) with \eqref{eq:M-decay} and
\(q(r)\le Mr^{1-N}\), we conclude that
\begin{equation}\label{eq:compactbound}
 M\le C,\qquad
 Q(r)\le Cr^{2-N},\qquad
 q(r)\le Cr^{1-N}\quad(r\ge1),
\end{equation}
where \(C>0\) depends only on \(N,p_0\).
\end{proof}

We next prove uniqueness and differentiable dependence on the value
at the origin.

\begin{lemma}\label{lem:radial-ivp}
Let \(N\ge3\).
For fixed \(p>2^*\) and \(\alpha>0\), the radial initial
value problem for \eqref{eq:pde}, with
\(u(0)=\alpha\) and \(u'(0)=0\), has a unique positive
local solution \(u_\alpha\). Uniqueness holds on every interval
where the solution remains positive. If \(u_\alpha>0\) on
\([0,R]\), then, for \(\gamma\) sufficiently close to \(\alpha\),
the solution \(u_\gamma\) is defined and positive on \([0,R]\),
and the map
\[
 \gamma\longmapsto u_\gamma|_{[0,R]}
\]
is \(C^1\) with values in \(C^1([0,R])\).
\end{lemma}

\begin{proof}
For a radial solution \(Q\), the flux formula \eqref{eq:flux1}
can be written as
\[
 m(r)=r^N\int_0^1 t^{N-1}Q(rt)^{p-1}\,dt.
\]
Substitution into \eqref{eq:compact-flux} gives
\begin{equation}\label{eq:ivp}
 Q'(r)=
 -\frac{r\int_0^1 t^{N-1}Q(rt)^{p-1}\,dt}
 {\sqrt{1+\left(r\int_0^1 t^{N-1}Q(rt)^{p-1}\,dt\right)^2}}.
\end{equation}
For an initial value \(\gamma\) near \(\alpha\), define
\[
 \mathcal V_\gamma(v)(r)
 =\gamma-\int_0^r
 \frac{s\int_0^1 t^{N-1}v(st)^{p-1}\,dt}
 {\sqrt{1+\left(s\int_0^1 t^{N-1}v(st)^{p-1}\,dt\right)^2}}\,ds.
\]
Consider the closed subset
\[
 \mathcal K_d=\{v\in C([0,d]):
                 \alpha/2\le v(r)\le2\alpha\text{ for }0\le r\le d\}.
\]
Since \(|\mathcal V_\gamma(v)(r)-\gamma|\le r\), the map
\(\mathcal V_\gamma\) takes \(\mathcal K_d\) into itself
when \(|\gamma-\alpha|\le\alpha/4\) and \(d\le\alpha/4\).
The function \(z\mapsto z/\sqrt{1+z^2}\) is
\(1\)-Lipschitz. Also,
\[
 |a^{p-1}-b^{p-1}|
 \le(p-1)(2\alpha)^{p-2}|a-b|
 \qquad(a,b\in[\alpha/2,2\alpha]).
\]
Hence, for \(v,w\in\mathcal K_d\),
\begin{equation}\label{eq:ivp-local}
 \begin{aligned}
 \|\mathcal V_\gamma(v)-\mathcal V_\gamma(w)\|_\infty
 &\le\frac{(p-1)(2\alpha)^{p-2}}{N}
       \int_0^d s\,ds\,\|v-w\|_\infty\\
 &=\frac{(p-1)(2\alpha)^{p-2}}{2N}
       d^2\|v-w\|_\infty.
 \end{aligned}
\end{equation}
Choose \(d\le\alpha/4\) so that
\((p-1)(2\alpha)^{p-2}d^2/(2N)<1/2\).
The choice of \(d\) depends only on \(N,p,\alpha\), not on
\(\gamma\). The contraction theorem gives a unique fixed point
\(u_\gamma\in\mathcal K_d\).

The map \((\gamma,v)\mapsto v-\mathcal V_\gamma(v)\) is
\(C^1\) near \((\alpha,u_\alpha)\).
By \eqref{eq:ivp-local}, the operator
\(I-D_v\mathcal V_\alpha(u_\alpha)\) is invertible.
The implicit function theorem therefore gives \(C^1\)
dependence of the solution on \(\gamma\) in \(C([0,d])\).
The right-hand side of \eqref{eq:ivp} is \(C^1\) in the positive
profile, so \(\gamma\mapsto u_\gamma'\) is also \(C^1\) in
\(C([0,d])\). Hence \(\gamma\mapsto u_\gamma\) is \(C^1\)
with values in \(C^1([0,d])\).
For \(r>0\), the equations
\[
 u'(r)=-\frac{m(r)}{\sqrt{r^{2N-2}+m(r)^2}},
 \qquad m'(r)=r^{N-1}u(r)^{p-1}
\]
are regular as long as \(u>0\). Uniqueness and \(C^1\)
dependence on \(\gamma\) thus extend to every compact
interval \([0,R]\) on which \(u_\alpha\) is positive.
For \(\gamma\) sufficiently close to \(\alpha\),
\(u_\gamma\) remains positive on \([0,R]\).
\end{proof}

We use the uniform bounds to pass to a limit when the exponent varies.

\begin{lemma}\label{lem:radial-compactness}
Let \(N\ge3\), \(p_i\to\bar p>2^*\), and let \(Q_i\) be positive radial
finite-energy solutions of \eqref{eq:pde} with exponents \(p_i\).
There are a subsequence, still denoted by \(Q_i\), and a positive
radial finite-energy solution \(Q\) with exponent \(\bar p\)
such that
\[
 \|Q_i-Q\|_{C^1([0,R])}\longrightarrow0
 \qquad\text{for every }R>0.
\]
\end{lemma}

\begin{proof}
Choose \(p_0\in(2^*,\bar p)\). After discarding finitely
many terms, we have \(p_i\ge p_0\).
By Lemma~\ref{lem:radial-bounds}, the values \(Q_i(0)\)
are bounded above and bounded away from zero. Since \(|Q_i'|\le1\),
a subsequence converges uniformly on every
bounded interval to a function \(Q\ge0\) with \(Q(0)>0\).
The same subsequence satisfies
\[
 Q_i^{p_i-1}\longrightarrow Q^{\bar p-1}
 \quad\text{uniformly on every bounded interval}.
\]
The convergence of the powers follows from local uniform
convergence of \(Q_i\), \(p_i\to\bar p\), and continuity of
\((p,z)\mapsto z^{p-1}\) for \(p>2\), \(z\ge0\).
Equation \eqref{eq:ivp} now gives \(Q_i'\to Q'\) uniformly
on every bounded interval, including at zero. The limit satisfies
\[
 m(r)=\int_0^r t^{N-1}Q(t)^{\bar p-1}\,dt,
 \qquad
 Q'(r)=-\frac{m(r)}{\sqrt{r^{2N-2}+m(r)^2}}.
\]
If \(Q\) had a first zero \(r_*>0\), then \(m(r_*)>0\)
would imply \(Q'(r_*)<0\), contrary to \(Q\ge0\).
Thus \(Q(r)>0\) for every \(r\ge0\).

Passing to the limit in the decay bounds of
Lemma~\ref{lem:radial-bounds}, we obtain
\[
 Q(r)\le Cr^{2-N},\qquad |Q'(r)|\le Cr^{1-N}
 \qquad(r\ge1),
\]
with \(C\) depending only on \(N,p_0\).
The decay estimates give \(Q\in L^{2^*}(\mathbb R^N)\) and
\(\nabla Q\in L^2(\mathbb R^N)\).
The formula for \(Q'\) gives \(|Q'(r)|<1\) for \(r>0\),
and \eqref{eq:ivp} gives \(Q'(0)=0\).
Since \(Q'(r)\to0\) as \(r\to\infty\), we also have
\(\|\nabla Q\|_\infty<1\).
Finally, \(m(r)=O(r^N)\) near zero, so the radial equation
gives the weak formulation \eqref{eq:weak} on all of
\(\mathbb R^N\), with exponent \(\bar p\).
Therefore \(Q\) is a positive radial finite-energy solution.
\end{proof}

\subsection{The number of radial solutions of index one}
We count solutions by their values at the origin:
\begin{equation}\label{eq:count-set}
 \mathcal S_p=\{Q(0): Q\text{ is a positive radial finite-energy solution of }
       \eqref{eq:pde},\ \ind(Q)=1\}.
\end{equation}
Let \(\#\mathcal S_p\) denote the number of elements of \(\mathcal S_p\).
By Lemma~\ref{lem:radial-ivp}, a central value determines at most
one positive radial solution, so \(\#\mathcal S_p\) counts the
positive radial finite-energy solutions with \(\ind(Q)=1\).
A positive radial ground state exists by
\cite[Theorem 1.2(iii)]{bim} and has radial Morse index one by
Lemma~\ref{lem:ground-morse}. Hence \(\mathcal S_p\ne\varnothing\).

We now count solutions as \(p\) varies.

\begin{lemma}\label{lem:count}
The number \(\#\mathcal S_p\) is finite and independent of \(p\)
for \(p>2^*\).
\end{lemma}
\begin{proof}
Fix \(\bar p>2^*\). We prove that \(\mathcal S_p\) is finite
and that \(\#\mathcal S_p\) is constant for \(p\) near \(\bar p\).
Choose \(2^*<p_0<\bar p<p_1\), and restrict \(p\) to
\([p_0,p_1]\).

\medskip\noindent
Step 1. We write the radial equation as a fixed point equation
and prove continuity, differentiability, and compactness estimates
for the resulting map.

Choose \(N/(2^*-1)<\nu<N-2\), and define the Banach space
\[
 X=\left\{u\in C([0,\infty)):
 \|u\|_X:=\sup_{r\ge0}(1+r)^\nu|u(r)|<\infty\right\}.
\]
Since \((N-2)(2^*-1)=N+2>N\), such a \(\nu\) exists and
satisfies \(\nu(p_0-1)>N\).

We first show that local uniform convergence and a uniform decay
bound imply convergence in \(X\).
Suppose \(u_i\to u\) locally uniformly and
\(|u_i(r)|+|u(r)|\le C r^{2-N}\) for \(r\ge1\), with \(C\)
independent of \(i\). Splitting the supremum at \(R\ge1\) gives
\begin{equation}\label{eq:X-conv}
 \|u_i-u\|_X
 \le(1+R)^\nu\sup_{0\le r\le R}|u_i(r)-u(r)|
     +C R^{\nu+2-N}.
\end{equation}
First let \(i\to\infty\) at fixed \(R\), then let \(R\to\infty\).
Since \(\nu<N-2\), \eqref{eq:X-conv} gives \(u_i\to u\) in \(X\).
Lemma~\ref{lem:radial-bounds} supplies the required uniform decay
for the solution sequences used below.

For \(u\in X\), set \(u_+=\max\{u,0\}\).
Integrating the flux formulas \eqref{eq:flux1} and
\eqref{eq:flux2} with \(m(0)=u(\infty)=0\) leads to the map
\begin{equation}\label{eq:fixed}
 \begin{aligned}
 m_u(t)&=\int_0^t s^{N-1}u_+(s)^{p-1}\,ds,\\
 T_pu(r)&=\int_r^\infty
 \frac{m_u(t)}{\sqrt{t^{2N-2}+m_u(t)^2}}\,dt.
 \end{aligned}
\end{equation}
The flux \(m_u\) also depends on \(p\). Using \(u_+\) defines
\(T_p\) even when \(u\) changes sign.
Every positive radial finite-energy solution satisfies \(T_pu=u\).
If \(T_pu=u\) and \(u(0)>0\), continuity gives \(u>0\) near zero.
Hence \(m_u(t)>0\) for every \(t>0\). Formula \eqref{eq:fixed}
then gives
\[
 u(r)>0,\qquad
 -u'(r)=\frac{m_u(r)}{\sqrt{r^{2N-2}+m_u(r)^2}}
 \qquad(r>0).
\]
We next prove the bounds needed to show that such a fixed point
has finite energy.

Let \(u\) range over a bounded subset of \(X\), and let
\(\|h\|_X\le1\). In the estimates for \(T_p\) below, the constants
\(C\) may depend on \(N,p_0,p_1,\nu\) and the bound for \(\|u\|_X\),
but not on \(p\in[p_0,p_1]\), \(h\), or the radial variables.
Since \(z\mapsto z_+^{p-1}\) is \(C^1\), the directional
derivative of \(m_u\) is
\[
 J_h(t)=(p-1)\int_0^t s^{N-1}u_+(s)^{p-2}h(s)\,ds.
\]
The bounds \(|u(s)|+|h(s)|\le C(1+s)^{-\nu}\) give
\begin{equation}\label{eq:mapbound}
 \begin{aligned}
 |m_u(t)|+|J_h(t)|
 &\le C\int_0^t s^{N-1}(1+s)^{-\nu(p_0-1)}\,ds\\
 &\le C\min\{t^N,1\}.
 \end{aligned}
\end{equation}
The last integral is bounded as \(t\to\infty\) because
\(\nu(p_0-1)>N\).
For \(t>0\),
\[
 \frac{d}{ds}\left(\frac{s}{\sqrt{t^{2N-2}+s^2}}\right)
 =\frac{t^{2N-2}}{(t^{2N-2}+s^2)^{3/2}}
 \le t^{1-N}.
\]
Consequently, \eqref{eq:mapbound} gives
\begin{equation}\label{eq:map-integrands}
 \frac{m_u(t)}{\sqrt{t^{2N-2}+m_u(t)^2}}
 +\frac{t^{2N-2}|J_h(t)|}
 {(t^{2N-2}+m_u(t)^2)^{3/2}}
 \le C\min\{t,t^{1-N}\}.
\end{equation}
The right-hand side is integrable on \((0,\infty)\).
Differentiating along \(u+\tau h\) under the integral sign gives
\begin{equation}\label{eq:dermap}
 [DT_p(u)h](r)
 =\int_r^\infty
 \frac{t^{2N-2}}{(t^{2N-2}+m_u(t)^2)^{3/2}}J_h(t)\,dt.
\end{equation}
After proving continuity of \(DT_p(u)\), we will show that
\eqref{eq:dermap} is the Fr\'echet derivative in \(X\).

Integrating \eqref{eq:map-integrands} and differentiating
\eqref{eq:fixed} and \eqref{eq:dermap} with respect to \(r\),
we obtain
\begin{equation}\label{eq:mapregular}
 \begin{aligned}
 |T_pu(r)|+|[DT_p(u)h](r)|
 &\le C(1+r)^{2-N},\\
 |(T_pu)'(r)|+|[DT_p(u)h]'(r)|
 &\le C\min\{r,r^{1-N}\}.
 \end{aligned}
\end{equation}
In particular, \(T_pu\) and \(DT_p(u)h\) belong to \(X\).
For \(w=T_pu\) or \(w=DT_p(u)h\), the decay bound in
\eqref{eq:mapregular} gives \(w\in L^{2^*}(\R^N)\), and
\[
 \int_0^\infty r^{N-1}|w'(r)|^2\,dr
 \le C\left(\int_0^1 r^{N+1}\,dr
             +\int_1^\infty r^{1-N}\,dr\right)<\infty.
\]
Thus \(T_pu\) and \(DT_p(u)h\) belong to \(\D_{\rm rad}\).

If \(T_pu=u\) and \(u(0)>0\), \eqref{eq:mapbound} gives
\[
 \sup_{r>0}\frac{m_u(r)}{r^{N-1}}<\infty.
\]
Differentiating \eqref{eq:fixed} then gives
\(\|u'\|_\infty<1\).
Since \(m_u(r)=O(r^N)\) near zero, the boundary term at zero
vanishes when we integrate the radial equation by parts.
Hence \(u\) satisfies \eqref{eq:weak} on \(\R^N\).
Since \(u\in\D_{\rm rad}\) and \(\|u'\|_\infty<1\),
the fixed point \(u\) is a positive radial finite-energy solution.

We next prove continuity in \(p,u\).
Suppose \(p_i\to p\) in \([p_0,p_1]\) and \(u_i\to u\) in \(X\).
In the notation \(m_{u_i}\), use exponent \(p_i\).
For each \(t>0\), local uniform convergence gives
\(m_{u_i}(t)\to m_u(t)\).
For the inner integral in \eqref{eq:dermap},
\[
 \begin{aligned}
 \sup_{\|h\|_X\le1}\bigg|
 &\int_0^t s^{N-1}
 \big((p_i-1)(u_i)_+^{p_i-2}-(p-1)u_+^{p-2}\big)h\,ds
 \bigg|\\
 &\le\frac{t^N}{N}\sup_{0\le s\le t}
 \big|(p_i-1)(u_i(s))_+^{p_i-2}
       -(p-1)u_+(s)^{p-2}\big|\longrightarrow0.
 \end{aligned}
\]
The limit follows from continuity of \((p,z)\mapsto z_+^{p-2}\)
for \(p>2\), together with \(p_i\to p\) and local uniform
convergence of \(u_i\).
For each \(t>0\), the factor
\[
 \frac{t^{2N-2}}{(t^{2N-2}+m_{u_i}(t)^2)^{3/2}}
\]
also converges to the same expression with \(m_u(t)\).
Thus the integrands in \eqref{eq:dermap} converge uniformly
for \(\|h\|_X\le1\) at each fixed \(t>0\).
The bound \eqref{eq:map-integrands} and dominated convergence give
\[
 \sup_{\|h\|_X\le1}\sup_{0\le r\le R}
 \big|[(DT_{p_i}(u_i)-DT_p(u))h](r)\big|
 \longrightarrow0
 \qquad(R\ge1).
\]
On \([0,R]\), the weight \((1+r)^\nu\) is at most \((1+R)^\nu\).
On \([R,\infty)\), \eqref{eq:mapregular} gives the weighted bound
\(CR^{\nu+2-N}\) for both operators, uniformly for \(\|h\|_X\le1\).
Therefore,
\begin{equation}\label{eq:DT-cont}
 \limsup_{i\to\infty}
 \|DT_{p_i}(u_i)-DT_p(u)\|_{X\to X}
 \le CR^{\nu+2-N}\longrightarrow0
 \qquad(R\to\infty).
\end{equation}

The integrands in \eqref{eq:fixed} also converge pointwise and
satisfy \eqref{eq:map-integrands}. Dominated convergence gives
\(T_{p_i}u_i\to T_pu\) locally uniformly. The uniform decay in
\eqref{eq:mapregular} and \eqref{eq:X-conv} then give
\begin{equation}\label{eq:T-cont}
 \|T_{p_i}u_i-T_pu\|_X\longrightarrow0.
\end{equation}
For fixed \(p\), integrating the directional derivative along
\(u+\tau h\) and using \eqref{eq:DT-cont}, we obtain
\[
 \begin{aligned}
 &\|T_p(u+h)-T_pu-DT_p(u)h\|_X\\
 &\qquad\le\|h\|_X\sup_{0\le\tau\le1}
 \|DT_p(u+\tau h)-DT_p(u)\|_{X\to X}
 =o(\|h\|_X).
 \end{aligned}
\]
Hence \(T_p\) is Fr\'echet differentiable in \(u\), and
\eqref{eq:DT-cont} proves that \(DT_p(u)\) is continuous in \(p,u\).

For fixed \(p,u\), the operator \(DT_p(u):X\to X\) is compact.
Indeed, if \(\|h_i\|_X\le1\), \eqref{eq:mapregular} gives uniform
bounds and equicontinuity of \(DT_p(u)h_i\) on every bounded interval.
The Arzel\`a--Ascoli theorem gives a locally uniformly convergent
subsequence. The decay bound in \eqref{eq:mapregular} and
\eqref{eq:X-conv} give convergence of that subsequence in \(X\).

\medskip\noindent
Step 2. We show that each solution with radial Morse index one
has a unique nearby solution when the exponent changes slightly.

Let \(Q\) be a solution with \(Q(0)\in\mathcal S_{\bar p}\),
and set \(q=-Q'\). Suppose \(h\in X\) satisfies
\(h=DT_{\bar p}(Q)h\).
By \eqref{eq:mapregular} and linearity in \(h\),
we have \(h\in\D_{\rm rad}\).
The flux formula \eqref{eq:flux2} gives
\[
 \frac{r^{2N-2}}{(r^{2N-2}+m_Q(r)^2)^{3/2}}
 =r^{1-N}(1-q(r)^2)^{3/2}.
\]
Differentiating \eqref{eq:dermap} with respect to \(r\), we obtain
\[
 \begin{aligned}
 -r^{N-1}(1-q^2)^{-3/2}h'&=J_h,\\
 J_h'&=(\bar p-1)r^{N-1}Q^{\bar p-2}h,
 \qquad J_h(0)=0.
 \end{aligned}
\]
The condition \(J_h(0)=0\) removes the boundary term at the origin.
For radial \(Q,h\), the bilinear form \eqref{eq:form} has the same
value on a test function and its spherical average.
Thus the equations for \(h,J_h\) imply \(L_Qh=0\) on \(\R^N\).
Since \(\ind(Q)=1\), Lemma~\ref{lem:rad} gives \(h=0\).
The operator \(DT_{\bar p}(Q)\) is compact by Step 1.
The Fredholm alternative therefore gives a bounded inverse for
\[
 A_Q:=I-DT_{\bar p}(Q).
\]

Define
\begin{equation}\label{eq:continuation-map}
 G_p(u):=u-A_Q^{-1}(u-T_pu).
\end{equation}
Then \(G_{\bar p}(Q)=Q\) and \(D_uG_{\bar p}(Q)=0\).
By \eqref{eq:DT-cont}, we can choose
\(0<\rho<Q(0)/2\) and a neighborhood of \(\bar p\) such that
\[
 \|D_uG_p(u)\|_{X\to X}\le\frac12
 \qquad\text{whenever }\|u-Q\|_X\le\rho.
\]
By \eqref{eq:T-cont}, \(T_pQ\to T_{\bar p}Q=Q\) in \(X\).
After reducing the neighborhood of \(\bar p\), we also have
\(\|G_p(Q)-Q\|_X\le\rho/2\). Therefore,
\[
 \|G_p(u)-Q\|_X
 \le\frac12\|u-Q\|_X+\frac\rho2
 \le\rho
 \qquad\text{if }\|u-Q\|_X\le\rho.
\]
Thus \(G_p\) is a contraction of the closed ball
\(\{u:\|u-Q\|_X\le\rho\}\) into itself.
Since \(A_Q\) is invertible, \(G_p(u)=u\) is equivalent to
\(T_pu=u\). The contraction theorem gives a unique \(u_p\) such that
\begin{equation}\label{eq:local-branch}
 \{u\in X:\|u-Q\|_X\le\rho,\ T_pu=u\}=\{u_p\}.
\end{equation}
Moreover, \(u_p(0)\ge Q(0)-\rho>0\).
By Step 1, \(u_p\) is a positive radial finite-energy solution.
The contraction estimate also gives
\begin{equation}\label{eq:branch-convergence}
 \|u_p-Q\|_X
 \le2\|A_Q^{-1}(Q-T_pQ)\|_X\longrightarrow0
 \qquad(p\to\bar p).
\end{equation}

\medskip\noindent
Step 3. We prove that the nearby solutions \(u_p\) have radial
Morse index one, and that limits of solutions with radial Morse
index one also have radial Morse index one.

We first prove convergence of the quadratic forms.
Let \(p_i\to p_*\) in \([p_0,p_1]\), and let \(Q_i\to Q_*\)
be a subsequence given by Lemma~\ref{lem:radial-compactness}.
Set \(q_i=-Q_i'\) and \(q_*=-Q_*'\).
For every \(R\ge1\), convergence in \(C^1([0,R])\) gives
uniform convergence of \(q_i\) on \([0,R]\).
Lemma~\ref{lem:radial-bounds} gives
\[
 |q_i(r)|+|q_*(r)|\le CR^{1-N}
 \qquad(r\ge R),
\]
where \(C\) is independent of \(i,R\).
Hence \(q_i\to q_*\) uniformly on \([0,\infty)\).
Since \(\|q_*\|_\infty<1\), we have
\[
 (1-q_i^2)^{-3/2}\longrightarrow(1-q_*^2)^{-3/2}
 \qquad\text{uniformly on }[0,\infty).
\]
The uniform decay bound for \(Q_i\) gives, for \(r\ge1\),
\[
 (p_i-1)Q_i(r)^{p_i-2}
 \le C r^{-(N-2)(p_0-2)},
 \qquad
 \int_1^\infty r^{N-1-(N-2)(p_0-2)N/2}\,dr<\infty.
\]
Here \(C\) is independent of \(i,r\), and
\((N-2)(p_0-2)>4\).
Local uniform convergence and dominated convergence therefore give
\[
 (p_i-1)Q_i^{p_i-2}
 \longrightarrow(p_*-1)Q_*^{p_*-2}
 \qquad\text{in }L^{N/2}(\R^N).
\]
By H\"older's inequality and the Sobolev inequality,
\begin{equation}\label{eq:formconv}
 \begin{aligned}
 &\sup_{\substack{h\in\D_{\rm rad}\\\|\nabla h\|_2=1}}
 |\B_{Q_i}(h,h)-\B_{Q_*}(h,h)|\\
 &\quad\le
 \|(1-q_i^2)^{-3/2}-(1-q_*^2)^{-3/2}\|_\infty\\
 &\qquad+C\|(p_i-1)Q_i^{p_i-2}
              -(p_*-1)Q_*^{p_*-2}\|_{N/2}
 \longrightarrow0.
 \end{aligned}
\end{equation}

The same convergence holds for \(Q_i=u_{p_i}\), \(Q_*=Q\),
and \(p_i\to\bar p\).
Indeed, \eqref{eq:branch-convergence} gives \(u_{p_i}\to Q\)
in \(X\). Equation \eqref{eq:ivp} gives local uniform convergence
of the derivatives. Lemma~\ref{lem:radial-bounds} bounds the
derivatives for large \(r\). Thus the coefficient and potential
estimates used in \eqref{eq:formconv} also apply to \(u_{p_i}\).

Return to the fixed solution \(Q\) at exponent \(\bar p\)
from Step 2. On \(\D_{\rm rad}\), define
\[
 (h,w)_Q
 =\int_{\R^N}(1-q^2)^{-3/2}\nabla h\cdot\nabla w\dx,
 \qquad
 \|h\|_Q^2=(h,h)_Q.
\]
Since \(\|q\|_\infty<1\), the norm \(\|\cdot\|_Q\) is
equivalent to the Dirichlet norm.
Define the operator \(K\) on \(\D_{\rm rad}\) by
\[
 (Kh,w)_Q=(\bar p-1)\int_{\R^N}Q^{\bar p-2}hw\dx.
\]
Then \(K\) is self-adjoint and nonnegative, and
\[
 \B_Q(h,w)=((I-K)h,w)_Q.
\]
We show that \(K\) is compact. For fixed \(R>0\), define
\[
 (K_Rh,w)_Q=(\bar p-1)\int_{B_R}Q^{\bar p-2}hw\dx.
\]
The Sobolev inequality bounds every bounded subset of
\(\D(\R^N)\) in \(H^1(B_R)\). The Rellich theorem therefore
makes the restriction map \(\D(\R^N)\to L^2(B_R)\) compact.
Since \(Q^{\bar p-2}\) is bounded on \(B_R\), the operator
\(K_R\) is compact.
On the complement of \(B_R\),
H\"older's inequality and the Sobolev inequality give
\[
 \begin{aligned}
 &(\bar p-1)\left|\int_{|x|>R}Q^{\bar p-2}hw\dx\right|\\
 &\qquad\le C(\bar p-1)
 \|Q^{\bar p-2}\|_{L^{N/2}(|x|>R)}\|h\|_Q\|w\|_Q.
 \end{aligned}
\]
Since \(Q^{\bar p-2}\in L^{N/2}(\R^N)\), the coefficient
on the right tends to zero as \(R\to\infty\).
Hence \(K_R\to K\) in operator norm, and \(K\) is compact.

Lemma~\ref{lem:rad} implies that \(1\) is not an eigenvalue of \(K\).
Since \(\ind(Q)=1\), exactly one eigenvalue of \(K\), counted
with multiplicity, is greater than \(1\).
The spectral theorem gives a vector \(e\) with \(\|e\|_Q=1\)
and a number \(\delta>0\) such that
\begin{equation}\label{eq:spectral-gap}
 \B_Q(e,e)\le-\delta,
 \qquad
 \B_Q(h,h)\ge\delta\|h\|_Q^2
 \quad\text{if }(h,e)_Q=0.
\end{equation}
Since \(\|\nabla h\|_2\le\|h\|_Q\), \eqref{eq:formconv}
gives, for \(p\) sufficiently close to \(\bar p\),
\[
 |\B_{u_p}(h,h)-\B_Q(h,h)|
 \le\frac\delta2\|h\|_Q^2
 \qquad(h\in\D_{\rm rad}).
\]
Combining with \eqref{eq:spectral-gap}, we obtain
\[
 \B_{u_p}(e,e)\le-\frac\delta2,
 \qquad
 \B_{u_p}(h,h)\ge\frac\delta2\|h\|_Q^2
 \quad\text{if }(h,e)_Q=0.
\]
The first inequality gives \(\ind(u_p)\ge1\).
Every two-dimensional subspace of \(\D_{\rm rad}\) contains
a nonzero \(h\) with \((h,e)_Q=0\), so the second inequality
gives \(\ind(u_p)\le1\). Therefore \(\ind(u_p)=1\).

We also need to preserve the index when taking limits.
Suppose \(\ind(Q_i)=1\) and \(Q_i\to Q_*\) as in
Lemma~\ref{lem:radial-compactness}.
If \(\B_{Q_*}\) were negative definite on a two-dimensional
space \(V\subset\D_{\rm rad}\), compactness of the unit circle
in \(V\) would give
\[
 \max_{\substack{h\in V\\\|\nabla h\|_2=1}}
 \B_{Q_*}(h,h)<0.
\]
By \eqref{eq:formconv}, \(\B_{Q_i}\) would then be negative
definite on \(V\) for all sufficiently large \(i\), contrary
to \(\ind(Q_i)=1\). Thus \(\ind(Q_*)\le1\), while
Lemma~\ref{lem:negative-direction} gives \(\ind(Q_*)\ge1\).
Hence \(\ind(Q_*)=1\).

\medskip\noindent
Step 4. We prove that there are finitely many solutions at
\(\bar p\) and that their number is unchanged for nearby exponents.

Suppose there were infinitely many distinct solutions
\(Q_i\) with \(Q_i(0)\in\mathcal S_{\bar p}\).
Lemma~\ref{lem:radial-compactness} gives a subsequence converging
locally in \(C^1\) to a positive radial finite-energy solution \(Q_*\).
Step 3 gives \(\ind(Q_*)=1\), and \eqref{eq:X-conv} gives
\(Q_i\to Q_*\) in \(X\).
Applying \eqref{eq:local-branch} at \(Q_*\) with \(p=\bar p\),
we find that \(Q_*\) is the only solution in a sufficiently
small \(X\)-ball around \(Q_*\).
Thus \(Q_i=Q_*\) for all sufficiently large \(i\), contrary
to the choice of distinct solutions.
Therefore \(\mathcal S_{\bar p}\) is finite.

List the solutions with radial Morse index one at \(\bar p\)
as \(Q^1,\ldots,Q^m\).
Choose pairwise disjoint closed \(X\)-balls centered at
\(Q^1,\ldots,Q^m\) as in \eqref{eq:local-branch}.
By Steps 2 and 3, for every \(p\) sufficiently close to
\(\bar p\), each ball contains exactly one positive radial
finite-energy solution with radial Morse index one.

There are no additional solutions with radial Morse index one
outside these balls when \(p\) is close to \(\bar p\).
Otherwise, we could choose \(p_i\to\bar p\) and solutions
\(\widetilde Q_i\) outside all the balls.
Lemma~\ref{lem:radial-compactness} and Step 3 would give a
subsequence converging to a solution with radial Morse index one
at \(\bar p\), which must be one of \(Q^1,\ldots,Q^m\).
Equation \eqref{eq:X-conv} would give convergence in \(X\),
so \(\widetilde Q_i\) would eventually lie in the ball around
that solution, a contradiction.

Hence \(\#\mathcal S_p=m=\#\mathcal S_{\bar p}\) for \(p\)
near \(\bar p\).
Since \(\bar p>2^*\) was arbitrary, \(\#\mathcal S_p\) is
finite and locally constant on \((2^*,\infty)\).
The interval \((2^*,\infty)\) is connected, so
\(\#\mathcal S_p\) is constant on the whole interval.
\end{proof}

\section{Proof of uniqueness and the equality cases}\label{sec:uniqueproof}
By Lemma~\ref{lem:count}, it suffices to prove \(\#\mathcal S_p=1\)
for large \(p\).

\subsection{Dependence on the value at the origin}
Fix \(p_0>2^*\). For \(p\ge p_0\), let \(Q\) be a positive radial
finite-energy solution and set \(\alpha=Q(0)\).
The constants \(c,C>0\) in this section depend only on \(N,p_0\),
unless a different dependence is stated.
Lemma~\ref{lem:radial-bounds} gives \(\alpha\ge c>0\).
Using \eqref{eq:flux2} and \eqref{eq:gint}, we obtain
\[
 c\le\alpha=\int_0^\infty \frac{m(r)}{\sqrt{r^{2N-2}+m(r)^2}}\,dr
 \le\int_0^\infty \frac{M}{\sqrt{r^{2N-2}+M^2}}\,dr\le C M^{1/(N-1)} .
\]
Thus \(M\ge c>0\). Choose the radius \(R\) from
Lemma~\ref{lem:radial-bounds}, so that \(m(R)\ge3M/4\). If \(M/2\le m(r)\le3M/4\), then \(r\le R\) and
\[
 q(r)\ge \frac{M/2}{\sqrt{R^{2N-2}+(M/2)^2}}\ge c>0,\qquad
 M/2\le m(r)\le \alpha^{p-1}r^N/N .
\]
Hence \(rq\ge c\,\alpha^{-(p-1)/N}\) where
\(M/2\le m(r)\le3M/4\). The increase of \(m\) over this
interval is \(M/4\ge c\). By \(dm=r^{N-1}Q^{p-1}\,dr\), we have
\begin{equation}\label{eq:median}
 \begin{aligned}
 c\,\alpha^{-(p-1)/N}
 &\le\int_{\{M/2\le m(r)\le3M/4\}} rq\,dm
   \le\int_0^\infty rq\,dm\\
 &=-\frac1p\int_0^\infty r^N(Q^p)'\,dr
   =\frac Np\int_0^\infty r^{N-1}Q^p\,dr\le C/p .
 \end{aligned}
\end{equation}
Here \(r^NQ^p\to0\) at both endpoints by regularity and the decay
\(Q(r)=O(r^{2-N})\) in \eqref{eq:Q-decay}.
The last inequality follows
from the uniform \(L^p\) bound in Lemma~\ref{lem:radial-bounds}.
Equation \eqref{eq:median} gives \(\alpha^{p-1}\ge cp^N\).
Since \(\alpha\ge c\), we obtain \(\alpha^p\ge cp^N\).
Combining with Lemma~\ref{lem:radial-bounds}, choose \(C_*>1\),
depending only on \(N,p_0\), such that
\begin{equation}\label{eq:heightinterval}
 C_*^{-1}p^N\le\alpha^p\le C_*p^N .
\end{equation}
Fix \(r_0=1/8\).

We compare the values and fluxes at \(r_0\) as the central value increases.

\begin{lemma}\label{lem:inner}
Let \(u\) solve \eqref{eq:ivp} with \(u(0)=\alpha\), where
\(\alpha\) satisfies \eqref{eq:heightinterval}. Let \(m\) be
its radial flux, defined by \eqref{eq:radial-flux-def} with
\(u\) in place of \(Q\).
For all sufficiently large \(p\), \(u>0\) on \([0,r_0]\) and
\begin{equation}\label{eq:signs}
 u(r_0)\le\theta<1,\qquad
 \partial_\alpha m(r_0)\ge cp,\qquad
 \partial_\alpha u(r_0)\le-cp .
\end{equation}
The constants \(c>0\) and \(\theta<1\) do not depend on
\(p\) or \(\alpha\) in this range.
\end{lemma}
\begin{proof}
By \eqref{eq:heightinterval},
\[
 |p\log\alpha-N\log p|\le\log C_* .
\]
Thus \(\alpha\to1\) uniformly for central values satisfying
\eqref{eq:heightinterval}. While \(u>0\),
we have \(-u'=m/\sqrt{r^{2N-2}+m^2}\le1\). For large \(p\)
and \(r\le r_0\), this gives
\(1/2<\alpha-r\le u(r)\le\alpha\) and
\(m(r)\le r_0^N\alpha^{p-1}/N\). Thus \(u\) stays positive
and \(m\) stays finite, so the initial value solution extends
throughout \([0,r_0]\).
As in \eqref{eq:core-height}, the bound \(u(r)\ge\alpha-r\)
gives
\[
 m(\alpha/p)\ge\frac{\alpha^{p+N-1}}{eNp^N}
 \ge c'>0,
\]
where the last inequality follows from \eqref{eq:heightinterval}
and \(\alpha\to1\).
For \(\alpha/p\le r\le r_0\), this gives
\[
 \frac{m(r)}{\sqrt{r^{2N-2}+m(r)^2}}
 \ge\frac{c'}{\sqrt{r_0^{2N-2}+(c')^2}}>0.
\]
Integrating this lower bound for \(-u'\) gives
\begin{equation}\label{eq:height-slope}
 u(r)\le\alpha-c(r-\alpha/p)_+ .
\end{equation}
For large \(p\), take
\(\alpha/p\le r_0/2\) and \(\alpha\le1+cr_0/4\). Then
\[
 u(r_0)\le1-cr_0/4=:\theta<1 .
\]
For \(0\le r\le r_0\), \eqref{eq:height-slope} and \(1-z\le e^{-z}\)
give
\begin{equation}\label{eq:envelope}
 \begin{aligned}
 u(r)^{p-1}&\le\alpha^{p-1}
  \exp[-c(p-1)(r-\alpha/p)_+/\alpha]\le Cp^N e^{-c'pr},\\
 u(r)^{p-2}&\le Cp^N e^{-c'pr}.
 \end{aligned}
\end{equation}
Here \(\alpha^{p-1}\le Cp^N\) by \eqref{eq:heightinterval}, and
\(1/2\le u\le\alpha\le2\) gives
\[
 \frac{p-1}{\alpha}(r-\alpha/p)_+\ge \frac{pr}{4}-1,\qquad
 u^{p-2}\le2u^{p-1}.
\]
The integral formula \eqref{eq:flux1} and the bound
\eqref{eq:envelope} give
\begin{equation}\label{eq:innerfluxbound}
 c'\le m(r)\le Cp^N\int_0^{r_0}t^{N-1}e^{-c'pt}\,dt
 \le C\int_0^\infty s^{N-1}e^{-c's}\,ds
 \qquad(\alpha/p\le r\le r_0).
\end{equation}

By Lemma~\ref{lem:radial-ivp}, we may set
\(h=\partial_\alpha u\) and \(j=\partial_\alpha m\).
Differentiating \eqref{eq:flux1} and \eqref{eq:flux2} with respect
to \(\alpha\) gives
\begin{equation}\label{eq:derinner}
 h(0)=1,\quad j(0)=0,\qquad
 h'=-\frac{r^{2N-2}}{(r^{2N-2}+m^2)^{3/2}}j,\qquad dj=(p-1)(h/u)\,dm .
\end{equation}
We first show that \(h\) stays close to one near the origin.
Let \(\varepsilon=p^{-1/2}\). For large \(p\),
\(\alpha/p<\varepsilon<r_0/2\). Since \(u\ge1/2\),
\eqref{eq:derinner} gives
\[
 |j(r)|\le(p-1)\int_0^r\frac{|h|}{u}\,dm
 \le 2p\,m(r)\sup_{[0,\varepsilon]}|h|\qquad(r\le \varepsilon).
\]
For \(m>0\), we have
\[
 \frac{m r^{2N-2}}{(r^{2N-2}+m^2)^{3/2}}=\frac{m/r^{N-1}}{(1+(m/r^{N-1})^2)^{3/2}}
 \le\min\{1,(m/r^{N-1})^{-2}\}.
\]
On \(0<r\le\alpha/p\) we have
\[
 \frac{m}{r^{N-1}}\ge \frac{r}{N}[\alpha(1-1/p)]^{p-1}
                \ge c r\alpha^p .
\]
On \([\alpha/p,\varepsilon]\), \eqref{eq:innerfluxbound} gives \(m\ge c\).
Consequently,
\begin{equation}\label{eq:core-g}
 \frac{m r^{2N-2}}{(r^{2N-2}+m^2)^{3/2}}\le
 \begin{cases}
 \min\{1,(c r\alpha^p)^{-2}\}, & 0<r\le\alpha/p,\\
 Cr^{2N-2}, & \alpha/p\le r\le\varepsilon .
 \end{cases}
\end{equation}
Since \(m\ge0\), \eqref{eq:derinner} gives
\[
 |h'|\le2p\,\frac{m r^{2N-2}}{(r^{2N-2}+m^2)^{3/2}}\sup_{[0,\varepsilon]}|h|.
\]
By \eqref{eq:core-g} and
\(\int_0^\infty\min\{1,(c r\alpha^p)^{-2}\}\dr=2/(c\alpha^p)\),
\begin{equation}\label{eq:core}
 \begin{aligned}
 \sup_{[0,\varepsilon]}|h-1|
 &\le Cp\sup_{[0,\varepsilon]}|h|
 \left(\int_0^{\alpha/p}\min\{1,(c r\alpha^p)^{-2}\}\dr
           +C\int_{\alpha/p}^{\varepsilon}r^{2N-2}\dr\right)\\
 &\le Cp(\alpha^{-p}+p^{-(2N-1)/2})\sup_{[0,\varepsilon]}|h|.
 \end{aligned}
\end{equation}
By \eqref{eq:heightinterval}, the coefficient in \eqref{eq:core}
is at most \(C(p^{1-N}+p^{1-(2N-1)/2})\to0\), uniformly
in \(\alpha\). Using \(\sup|h|\le1+\sup|h-1|\) and
absorbing the term containing \(\sup|h-1|\), we obtain
\[
 \tfrac12\le h\le2\quad\text{on }[0,\varepsilon]
\]
for all sufficiently large \(p\).
By \eqref{eq:derinner}, \(1/2\le u\le2\), and \eqref{eq:innerfluxbound},
\begin{equation}\label{eq:coremass}
 \tfrac{p-1}{4}m(\varepsilon)\le j(\varepsilon)
 \le4(p-1)m(\varepsilon),\qquad cp\le j(\varepsilon)\le Cp .
\end{equation}

We next show that \(j\) stays positive on \([\varepsilon,r_0]\).
Equations \eqref{eq:envelope} and \eqref{eq:innerfluxbound} give
\[
 (p-1)u^{p-2}\le Cp^{N+1}e^{-c'\sqrt p}\le Cp^{-4},\qquad
 0<\frac{r^{2N-2}}{(r^{2N-2}+m^2)^{3/2}}\le Cr^{2N-2}.
\]
With all suprema taken over \([\varepsilon,r_0]\), \eqref{eq:derinner} gives
\begin{equation}\label{eq:outer-system}
 \begin{aligned}
 |j(r)-j(\varepsilon)|&\le
   \int_{\varepsilon}^r (p-1)t^{N-1}u^{p-2}|h|\,dt
   \le Cp^{-4}\sup|h|,\\
 |h(r)|&\le |h(\varepsilon)|+
      \int_{\varepsilon}^r \frac{t^{2N-2}}{(t^{2N-2}+m(t)^2)^{3/2}}|j(t)|\,dt
   \le 2+C\sup|j| .
 \end{aligned}
\end{equation}
By \eqref{eq:outer-system} and \eqref{eq:coremass},
\[
 (1-Cp^{-4})\sup|h|\le2+C|j(\varepsilon)|\le Cp,
 \qquad
 \sup|j-j(\varepsilon)|\le Cp^{-4}\sup|h|\le Cp^{-3}.
\]
Hence
\begin{equation}\label{eq:jpositive}
 \sup |h|\le Cp,\qquad j(r)\ge j(\varepsilon)-Cp^{-3}\ge cp
 \quad(\varepsilon\le r\le r_0).
\end{equation}
Then \(j(r_0)\ge cp\). For \(r_0/2\le r\le r_0\), we have \(m\le C\)
by \eqref{eq:innerfluxbound}, so
\[
 \frac{r^{2N-2}}{(r^{2N-2}+m^2)^{3/2}}\ge c>0 .
\]
Integrating the equation for \(h'\) in \eqref{eq:derinner}
and using \eqref{eq:jpositive}, we obtain
\[
 \begin{aligned}
 h(r_0)&=h(\varepsilon)-\int_{\varepsilon}^{r_0}\frac{r^{2N-2}}{(r^{2N-2}+m^2)^{3/2}}j\,dr\\
 &\le 2-\int_{r_0/2}^{r_0}\frac{r^{2N-2}}{(r^{2N-2}+m^2)^{3/2}}j\,dr
 \le 2-cp\,r_0/2 .
 \end{aligned}
\]
For large \(p\), \(h(r_0)\le-(cr_0/4)p\), which completes
the proof of \eqref{eq:signs}.
\end{proof}

\subsection{Comparison for large exponents}
\begin{proof}[Proof of Theorem~\ref{th:uniqueness}]
Let \(p\) be large. Suppose there are two values
\(\alpha_1<\alpha_2\) in \(\mathcal S_p\), and let
\(Q_i,m_i\) be the corresponding solutions and fluxes.
Every \(\alpha\in[\alpha_1,\alpha_2]\) satisfies
\eqref{eq:heightinterval}. Lemma~\ref{lem:inner} gives a positive
solution on \([0,r_0]\) for each such central value.
Only the solutions at the endpoints \(\alpha_1,\alpha_2\) need
to have finite energy on \(\R^N\); the intermediate initial
value solutions are used only on \([0,r_0]\).
Lemma~\ref{lem:radial-ivp} gives \(C^1\) dependence on
\(\alpha\). Integrating the bounds \eqref{eq:signs} over
\([\alpha_1,\alpha_2]\), we obtain
\begin{equation}\label{eq:opposite}
 \begin{aligned}
 m_2(r_0)-m_1(r_0)
 &=\int_{\alpha_1}^{\alpha_2}\partial_\alpha m(r_0;\alpha)\,d\alpha
 \ge cp(\alpha_2-\alpha_1)>0,\\
 Q_2(r_0)-Q_1(r_0)
 &=\int_{\alpha_1}^{\alpha_2}\partial_\alpha u(r_0;\alpha)\,d\alpha
 \le-cp(\alpha_2-\alpha_1)<0 .
 \end{aligned}
\end{equation}

Let \(v=(Q_1-Q_2)_+\) on \([r_0,\infty)\). Then \(v(r)\to0\) as
\(r\to\infty\), and almost everywhere
\[
 v'=(Q_1'-Q_2')\mathbf{1}_{\{Q_1>Q_2\}} .
\]
The Dirichlet bounds in Lemma~\ref{lem:radial-bounds} give
\(\int_{r_0}^\infty r^{N-1}|v'|^2\,dr<\infty\).
Using \eqref{eq:radial-flux-def} and \eqref{eq:flux1}, and integrating
by parts on \([r_0,S]\), we obtain
\begin{equation}\label{eq:exterior-finite}
 \int_{r_0}^S r^{N-1}(Q_1^{p-1}-Q_2^{p-1})v\,dr
 =[(m_1-m_2)v]_{r_0}^S
  +\int_{r_0}^S r^{N-1}\left(\frac{Q_1'}{\sqrt{1-(Q_1')^2}}-\frac{Q_2'}{\sqrt{1-(Q_2')^2}}\right)v'\,dr .
\end{equation}
Since each \(m_i\) has a finite limit by \eqref{eq:fast1} and
\(v(S)\to0\), the boundary term
\((m_1(S)-m_2(S))v(S)\) tends to zero.
Since \(\max_i\|Q_i'\|_\infty<1\), the derivative of
\(s/\sqrt{1-s^2}\) is bounded on the interval containing
the ranges of \(Q_1'\) and \(Q_2'\).
The bound may depend on the fixed pair \(Q_1,Q_2\).
The two integrands in
\eqref{eq:exterior-finite}, apart from the weight \(r^{N-1}\),
are bounded in absolute value by \(Q_1^p+Q_2^p\) and
\(C(|Q_1'|^2+|Q_2'|^2)\), respectively.
Lemma~\ref{lem:radial-bounds} shows that both bounds are
integrable with weight \(r^{N-1}\). Letting \(S\to\infty\), we have
\begin{equation}\label{eq:exteriorenergy}
 \begin{aligned}
 \int_{r_0}^\infty r^{N-1}\left(\frac{Q_1'}{\sqrt{1-(Q_1')^2}}-\frac{Q_2'}{\sqrt{1-(Q_2')^2}}\right)v'\,dr
 &= [m_1(r_0)-m_2(r_0)]v(r_0)\\
 &\quad+\int_{r_0}^\infty r^{N-1}(Q_1^{p-1}-Q_2^{p-1})v\,dr .
 \end{aligned}
\end{equation}
By \eqref{eq:opposite}, \(m_1(r_0)<m_2(r_0)\) and \(v(r_0)>0\), so
the boundary term is negative. On \(\{v>0\}\), by the mean value
theorem and \((1-s^2)^{-3/2}\ge1\), we have
\[
 \left(\frac{Q_1'}{\sqrt{1-(Q_1')^2}}-\frac{Q_2'}{\sqrt{1-(Q_2')^2}}\right)v'\ge |v'|^2,\qquad
 (Q_1^{p-1}-Q_2^{p-1})v\le(p-1)Q_1^{p-2}v^2.
\]
Both bounds also hold almost everywhere on \(\{v=0\}\).
By \eqref{eq:exteriorenergy}, we have
\begin{equation}\label{eq:energybound}
 \int_{r_0}^\infty r^{N-1}|v'|^2\dr
 \le(p-1)\int_{r_0}^\infty r^{N-1}Q_1^{p-2}v^2\dr .
\end{equation}
By \(v(\infty)=0\) and Cauchy--Schwarz,
\begin{equation}\label{eq:radialv}
 v(r)^2
 \le\left(\int_r^\infty t^{1-N}\,dt\right)
     \left(\int_r^\infty t^{N-1}|v'(t)|^2\,dt\right)
 \le\frac{r^{2-N}}{N-2}\int_{r_0}^\infty t^{N-1}|v'(t)|^2\,dt .
\end{equation}
By \eqref{eq:energybound} and \eqref{eq:radialv}, we have
\[
 \int_{r_0}^{\infty} r^{N-1}|v'|^2\dr
 \le\left(\frac{p-1}{N-2}\int_{r_0}^{\infty}rQ_1(r)^{p-2}\dr\right)
       \int_{r_0}^{\infty} r^{N-1}|v'|^2\dr .
\]
For \(r\ge r_0\), \(Q_1(r)\le\theta<1\) by monotonicity and
\eqref{eq:signs}. Also, \eqref{eq:M-decay} and the bound
\(M\le C\) in Lemma~\ref{lem:radial-bounds} give
\(Q_1(r)\le Cr^{2-N}\) for \(r\ge r_0\). Since \(N\ge3\), we have
\[
 \int_{r_0}^\infty r Q_1(r)^3\,dr
 \le C\int_{r_0}^\infty r^{1-3(N-2)}\,dr
 = C\,\frac{r_0^{8-3N}}{3N-8}<\infty
\]
uniformly in \(Q_1\). For \(p>5\), \(Q_1^{p-2}\le\theta^{p-5}Q_1^3\).
Then
\begin{equation}\label{eq:exterior-factor}
 \frac{p-1}{N-2}\int_{r_0}^\infty rQ_1(r)^{p-2}\,dr
 \le Cp\,\theta^{p-5}\int_{r_0}^\infty rQ_1(r)^3\,dr
 \le Cp\,\theta^{p-5}\longrightarrow0 .
\end{equation}
For large \(p\), the coefficient in \eqref{eq:exterior-factor}
is at most \(1/2\). Combining \eqref{eq:energybound},
\eqref{eq:radialv}, and \eqref{eq:exterior-factor} gives
\(\int_{r_0}^\infty r^{N-1}|v'|^2\dr=0\), and
\eqref{eq:radialv} gives \(v(r_0)=0\). This contradicts
\eqref{eq:opposite}.

Thus \(\#\mathcal S_p\le1\) for large \(p\), and existence of a
positive radial ground state gives \(\#\mathcal S_p=1\).
Lemma~\ref{lem:count} then gives \(\#\mathcal S_p=1\) for every
\(p>2^*\). By Lemmas~\ref{lem:radial-symmetry} and \ref{lem:ground-morse},
every positive ground state becomes radial with Morse index one
after translation to the origin. All such ground states then have
the unique central value in \(\mathcal S_p\), so
Lemma~\ref{lem:radial-ivp} makes them identical.
Strict decrease for \(r>0\) follows from \eqref{eq:radial-monotonicity}.
\end{proof}

\begin{proof}[Proof of Corollary~\ref{cor:eq}]
Let \(u\ne0\) attain equality in the sharp Born--Infeld Sobolev
inequality \eqref{eq:sharp}. By
\cite[Theorems 1.2(iii) and 1.4]{bim}, \(s^{-1}u(s\,\cdot)\)
is a ground state of constant sign for some \(s>0\). By
Theorem~\ref{th:uniqueness}, for some \(y\in\R^N\) and
\(\sigma\in\{-1,1\}\),
\[
 s^{-1}u(sx)=\sigma Q_p(x-y/s).
\]
Hence \(u\) has the form \eqref{eq:eqcases}.

Conversely, for every \(u\) in \eqref{eq:eqcases}, we have
\[
 \begin{aligned}
 \int_{\R^N}\bigl(1-\sqrt{1-|\nabla u|^2}\bigr)\dx
 &=s^N\int_{\R^N}\bigl(1-\sqrt{1-|\nabla Q_p|^2}\bigr)\dx,\\
 \|u\|_p^{Np/(N+p)}&=s^N\|Q_p\|_p^{Np/(N+p)} .
 \end{aligned}
\]
The function \(Q_p\) attains equality by \eqref{eq:ground-equality}.
Both sides of \eqref{eq:sharp} scale by \(s^N\), so \(u\)
also attains equality.
\end{proof}

\vspace{3mm}\par\noindent
{\bf Conflict of interest:} The author declares no conflict of interest.
\vspace{3mm}\par\noindent
{\bf Data availability statement:} Data sharing is not applicable to this article, as no datasets were generated or analyzed during this study.
\vspace{3mm}\par\noindent
{\bf
Declaration of generative AI}

During the preparation of this work, the author used GPT-6 Astra
to assist with checking proof steps and improving the English expression.
GPT-6 Astra also assisted with the proof of Lemma \ref{lem:rad} and Lemma \ref{lem:count}. The mathematical ideas were developed and verified by
the author. The author reviewed and revised the AI-generated material,
checked the proof arguments, and takes full responsibility for the
accuracy and rigor of the manuscript.

\end{document}